\documentclass[number,citesort,seceqn,dvips]{arxbj}
\usepackage{cursive}

\aid{0}
\volume{17}
\issue{1}
\pubyear{2011}
\firstpage{1}
\lastpage{33}
\doi{10.3150/10-BEJ259}

\makeatletter
\newcommand{\bbZ}{\mathbb{Z}}
\newcommand{\bbR}{\mathbb{R}}
\newcommand{\bbN}{\mathbb{N}}
\newcommand{\bbC}{\mathbb{C}}

\newcommand{\norm}[1]{\|#1\|}

\newproclaim{definition}{Definition}[section]
\newtheorem{theorem}{Theorem}[section]
\newremark{remark}{Remark}[section]
\newtheorem{proposition}{Proposition}[section]
\newremark{example}{Example}[section]
\newtheorem{corollary}{Corollary}[section]
\makeatother

\begin{document}
\begin{frontmatter}

\title{Integral representations and properties of operator fractional
Brownian motions}
\runtitle{Integral representations and properties of OFBMs}

\begin{aug}
\author[a]{\fnms{Gustavo} \snm{Didier}\thanksref{a}\ead[label=e1]{gdidier@tulane.edu}\corref{}}
\and
\author[b]{\fnms{Vladas} \snm{Pipiras}\thanksref{b}\ead[label=e2]{pipiras@email.unc.edu}}
\runauthor{G. Didier and V. Pipiras}
\address[a]{Mathematics Department,
Tulane University,
6823 St.\ Charles Avenue,
New Orleans, LA 70118, USA.
\printead{e1}}
\address[b]{Department\ of Statistics and Operations Research,
University of North Carolina at Chapel Hill,
CB\#3260, Smith Building,
Chapel Hill, NC 27599, USA.
\printead{e2}}
\end{aug}

\received{\smonth{9} \syear{2009}}
\revised{\smonth{1} \syear{2010}}

%
\begin{abstract}
Operator fractional Brownian motions (OFBMs) are (i) Gaussian,
(ii) operator self-similar and (iii) stationary increment
processes. They are the natural multivariate generalizations of
the well-studied fractional Brownian motions. Because of the
possible lack of time-reversibility, the defining properties
(i)--(iii) do not, in general, characterize the covariance
structure of OFBMs. To circumvent this problem, the class of OFBMs
is characterized here by means of their integral representations in
the spectral and time domains. For the spectral domain
representations, this involves showing how the operator
self-similarity shapes the spectral density in the general
representation of stationary increment processes. The time domain
representations are derived by using primary matrix functions and
taking the Fourier transforms of the deterministic spectral
domain kernels. Necessary and sufficient conditions for OFBMs to
be time-reversible are established in terms of their spectral and
time domain representations. It is also shown that the spectral
density of the stationary increments of an OFBM has a rigid
structure, here called the \textit{dichotomy principle}. The notion of operator
Brownian motions is also explored.
\end{abstract}

%

\end{frontmatter}

\section{Introduction}
\label{s:intro}

Fractional Brownian motion (FBM), denoted $B_{H} =\{B_H(t)\}_{t \in
\bbR}$ with $H \in(0,1)$, is a stochastic process characterized by
the following three properties:
\begin{longlist}[(iii)]
\item[(i)] Gaussianity;
\item[(ii)] self-similarity with parameter $H$;
\item[(iii)] stationarity of the increments.
\end{longlist}
By self-similarity, it is meant that the law of $B_H$ scales as
%
\begin{equation} \label{e:s.s.}
\{B_{H}(ct)\}_{t \in\bbR} \stackrel{{\mathcal L}}= \{c^{H}
B_{H}(t)\}_{t \in\bbR},\qquad  c > 0,
\end{equation}
where $\stackrel{{\mathcal L}}=$ denotes equality of finite-dimensional
distributions. By stationarity of the increments, it is meant that the
process
\[
\{B_{H}(t+h) - B_{H}(h)\}_{t \in\bbR}
\]
has the same distribution for any time-shift $h \in\bbR$. It may
be shown that these three properties actually characterize FBM in
the sense that it is the \textit{unique} (up to a constant) such
process for a given $H \in(0,1)$. FBM plays an important role in
both theory and applications, especially in connection with long-range
dependence \cite{embrechtsmaejima2002,doukhan2003}.

We are interested here in the multivariate counterparts of FBM,
called operator fractional Brownian motions (OFBMs). In the
multivariate context, an OFBM $B_H = (B_{1,H}, \ldots,
B_{n,H})^{*} = \{(B_{1,H}(t), \ldots, B_{n,H}(t))^{*} \in
\bbR^n, t \in\bbR\}$ is a collection of random vectors, where the
symbol $\ast$ denotes transposition. It is also Gaussian and has
stationary increments. Moreover, as is standard for the multivariate
context, in this paper, we assume that  OFBMs are \textit{proper},
that is, for each $t$, the distribution of $B_H(t)$ is not contained
in a proper subspace of $\bbR^n$. However, self-similarity is now
replaced by

(ii$'$) operator self-similarity.

A proper multivariate process $B_H$ is called (strictly) operator
self-similar (o.s.s.) if it is continuous in law for all $t$ and
the expression (\ref{e:s.s.}) holds for some matrix $H$. Here, the
expression $c^H$ is defined by means of the convergent series
\[
c^H = \exp(\log(c)H) = \sum_{k=0}^\infty(\log c)^k
\frac{H^k}{k!},\qquad  c> 0.
\]
Operator self-similarity extends the usual notion of self-similarity
and was first studied thoroughly in
\cite{laharohatgi1981,hudsonmason1982}; see also Section 11 in
\cite{meerschaertscheffler2001} and Chapter 9 in
\cite{embrechtsmaejima2002}. The theory of operator self-similarity
runs somewhat parallel to that of operator stable measures (see
\cite{jurekmason1993,meerschaertscheffler2001}) and is also related to
that of operator scaling random fields (see, e.g.,
\cite{biermemeerschaertscheffler2007}).

OFBMs are of interest in several areas and for reasons similar to
those in the univariate case. For example, OFBMs arise and are
used in the context of multivariate time series and long-range
dependence (see, e.g.,
\cite{marinuccirobinson2000,davidsondejong2000,chung2002,doladomarmol2004,davidsonhashimzade2007,robinson2008}).
Another context is that of queueing systems, where
reflected OFBMs model the size of multiple queues in particular
classes of queueing models and are studied in problems related
to, for example, large deviations (see \cite{konstantopouloslin1996,majewski2003,majewski2005,delgado2007}).
Partly motivated by this interest in OFBMs, several authors
consider constructions and properties of OFBMs. Maejima and Mason
\cite{maejimamason1994}, in particular, construct examples of
OFBMs through time domain integral representations. Mason and Xiao
\cite{masonxiao2002} study sample path properties of OFBMs.
Bahadoran, Benassi and D\c{e}bicki
\cite{bahadoranbenassidebicki2003} provide wavelet
decompositions of OFBMs, study their sample path properties and consider questions of identification. Becker-Kern and Pap
\cite{beckerkernpap2008} consider estimation of the real
spectrum of the self-similarity exponent. A number of other works
on operator self-similarity are naturally related to OFBMs; see,
for example,
\cite{meerschaertscheffler2001}, Section 11, and references therein.

To the reader less familiar with OFBMs, we should note that the
multivariate case is quite different from the univariate case. For
example, consider an OFBM $B_H$ whose exponent $H$ has
characteristic roots $h_k$ with positive real parts.
By using operator self-similarity and
stationarity of increments, one can argue, as in the univariate
case, that
%
\begin{eqnarray} \label{e:OFBM-cov}
&&EB_H(t)B_H(s)^{*} + EB_H(s)B_H(t)^{*}\nonumber
\\
&&\quad = EB_H(t)B_H(t)^{*}+EB_H(s)B_H(s)^{*} -
E\bigl(B_H(t)-B_H(s)\bigr)\bigl(B_H(t)-B_H(s)\bigr)^{*}\quad
\\
&&\quad = |t|^{H}\Gamma(1,1)|t|^{H^{*}}+
|s|^{H}\Gamma(1,1)|s|^{H^{*}}-|t-s|^{H}\Gamma(1,1)|t-s|^{H^{*}},\nonumber
\end{eqnarray}
where $\Gamma(1,1) = EB_H(1)B_H(1)^*$ and the symbol $*$ denotes
the adjoint operator. However, in contrast with the univariate case,
it is \textit{not} generally true that
%
\begin{equation}\label{e:EBH(t)BH(s)*_neq_EBH(s)BH(t)*}
EB_H(t)B_H(s)^{*} = EB_H(s)B_H(t)^{*}
\end{equation}
and hence the OFBM is not characterized by $H$ and a matrix
$\Gamma(1,1)$. Another important difference is that the
self-similarity exponent of an operator self-similar process is
generally not unique. The latter fact has been well known since
the fundamental work of Hudson and Mason \cite{hudsonmason1982}.
We briefly recall it, together with some related results, in Section
\ref{ss:o.s.s._processes}.

In this work, we address several new and, in our view, important
questions about OFBMs. In view of (\ref{e:OFBM-cov}) and
(\ref{e:EBH(t)BH(s)*_neq_EBH(s)BH(t)*}), since the covariance
structure of an OFBM cannot be determined in general, we pursue the
characterization of OFBMs in terms of their integral
representations (Section \ref{s:integral_repres_OFBM}). In the
spectral domain, under the mild and natural assumption that the
characteristic roots of $H$ satisfy
%
\begin{equation}\label{e:0<Re(hk)<1_intro}
0 < \operatorname{Re}(h_k) < 1,\qquad  k = 1,\dots,n,
\end{equation}
we show that an OFBM admits the integral
representation
%
\begin{equation}\label{e:OFBM-repres-sp}
\int_{\bbR} \frac{\mathrm{e}^{\mathrm{i}tx} - 1}{\mathrm{i}x} \bigl(x^{-(H-(1/2)I)}_{+}A +
x^{-(H-(1/2)I)}_{-}\overline{A}\bigr) \widetilde{B}(\mathrm{d}x),
\end{equation}
where $A$ is a matrix with complex-valued entries, $\overline{A}$
denotes its complex conjugate, $x_+ = \max\{x,0\}$, $x_- =
\max\{-x,0\}$ and $\widetilde B(\mathrm{d}x)$ is a suitable multivariate
complex-valued Gaussian measure. In the time domain and when, in
addition to (\ref{e:0<Re(hk)<1_intro}), we have
%
\begin{equation}\label{e:Re(hk)_neq_1/2_intro}
\operatorname{Re}(h_k) \neq\frac{1}{2},\qquad  k = 1,\dots,n,
\end{equation}
the OFBM admits the integral representation
%
\begin{equation}\label{e:OFBM-repres-time}
\hspace*{-3pt}\int_{\bbR} \bigl(\bigl((t - u)^{H-(1/2)I}_{+} - (-
u)^{H-(1/2)I}_{+}\bigr)M_{+} + \bigl((t - u)^{H-(1/2)I}_{-} - (-
u)^{H-(1/2)I}_{-}\bigr)M_{-} \bigr) B(\mathrm{d}u),\hspace*{3pt}
\end{equation}
where $M_+,M_-$ are matrices with real-valued entries and $B(\mathrm{d}u)$
is a suitable multivariate real-valued Gaussian measure. The
representation (\ref{e:OFBM-repres-time}) is obtained from
(\ref{e:OFBM-repres-sp}) by taking the Fourier transform of the
deterministic kernel in (\ref{e:OFBM-repres-sp}). We shall provide
rigorous arguments for this step by using primary matrix
functions. (Even in the univariate case, very often this step is
unjustifiably taken as more or less evident.) On a related note,
but from a different angle, the representations~(\ref{e:OFBM-repres-sp}) and (\ref{e:OFBM-repres-time}) always
define Gaussian processes with stationary increments that satisfy
(\ref{e:s.s.}) for a matrix $H$. We shall provide sufficient
condition for these processes to be proper (see Section~\ref{s:properness}) and, hence, to be OFBMs.

Subclasses of the representations (\ref{e:OFBM-repres-sp}) and
(\ref{e:OFBM-repres-time}) were considered in the works referenced
above. Maejima and Mason \cite{maejimamason1994} consider OFBMs
given by the representation (\ref{e:OFBM-repres-time}) with $M_+ =
M_- = I$. Mason and Xiao \cite{masonxiao2002} take
(\ref{e:OFBM-repres-sp}) with $A = I$. Bahadoran \textit{et al.}\
\cite{bahadoranbenassidebicki2003} consider
(\ref{e:OFBM-repres-sp}) with $A$ having full rank and real-valued
entries. (Such
OFBMs, for example, are necessarily time-reversible; see Theorem \ref
{t:time-revers_OFBM}
and also Remark \ref{r:example_proper_with_deficient_rank_A}.) We would
again like
to emphasize that, in contrast
with these works, the representations (\ref{e:OFBM-repres-sp}) and
(\ref{e:OFBM-repres-time}) characterize \textit{all} OFBMs (under
the mild and natural conditions (\ref{e:0<Re(hk)<1_intro}) and (\ref
{e:Re(hk)_neq_1/2_intro})).

In particular, the representations
(\ref{e:OFBM-repres-sp}) and (\ref{e:OFBM-repres-time}) provide a natural
framework for the study of many properties of OFBMs.
In this paper, we provide conditions
in terms of $A$ in (\ref{e:OFBM-repres-sp}) (or $M_+,M_-$ in
(\ref{e:OFBM-repres-time})) for OFBMs to be time-reversible (see
Section \ref{s:time_reversibility}). Time-reversibility is shown
to be equivalent to the condition (\ref
{e:EBH(t)BH(s)*_neq_EBH(s)BH(t)*}) and hence, in view of (\ref
{e:OFBM-cov}), corresponds to the situation
where the covariance structure of the OFBM is given by
%
\begin{equation}\label{e:OFBM-cov-time-revers}
EB_H(t)B_H(s)^* = \frac{1}{2} \bigl(|t|^{H}\Gamma(1,1)|t|^{H^{*}}+
|s|^{H}\Gamma(1,1)|s|^{H^{*}}-|t-s|^{H}\Gamma(1,1)|t-s|^{H^{*}} \bigr).
\end{equation}
Another interesting and little-explored direction of study of
OFBMs is their uniqueness (identification). This encompasses the
characterization of the different parameterizations for any given
OFBM and, in particular, of the aforementioned non-uniqueness of
the self-similarity exponents. Uniqueness questions in the
context of OFBMs are the focus of Didier and Pipiras
\cite{didierpipiras2010}, where they are explored starting with
the representation (\ref{e:OFBM-repres-sp}), and will be largely
absent from this paper.

Furthermore, in this paper, we also discuss some additional
properties of OFBMs which are of independent interest. First, we
prove that OFBMs have a rigid dependence structure among
components, which we call the \textit{dichotomy principle} (Section
\ref{s:dicho}). More precisely, under long-range dependence (in
the sense considered in Section \ref{s:dicho}), we show that the
components of the increments of an OFBM are either independent or
long-range dependent, that is, they cannot be short-range dependent
in a non-trivial way. Since, in the univariate case, the increments
of FBM are often considered representative of all long-range
dependent series, this result raises the question of whether OFBMs
are flexible enough to capture multivariate long-range dependence
structures. Second, we also discuss the notion of operator
Brownian motions (OBMs) and related questions (Section
\ref{s:OBM}). OBMs are defined as having independent increments
and are known to admit $H=(1/2)I$ as an exponent. We also show, in
particular, that an OFBM with $H=(1/2)I$ does not necessarily have
independent increments and hence is not necessarily an OBM. (In
contrast, in the univariate case, $H=1/2$ necessarily implies
Brownian motion.)

In summary, the structure of the paper is as follows. In Section \ref
{s:preliminaries},
we provide the necessary background for the paper and some definitions.
In Section \ref{s:integral_repres_OFBM}, we construct integral
representations for OFBMs
in the spectral and time domains. Section \ref{s:properness} furnishes
sufficient conditions for properness.
Section \ref{s:time_reversibility} is dedicated to time-reversibility.
The dichotomy principle is established in Section \ref{s:dicho} and the
properties
of OBMs are studied in Section \ref{s:OBM}. Appendices \ref{ss:FT}--\ref
{a:Jordan_form} contain several
important technical results used throughout the paper, as well as some proofs.

\section{Preliminaries}
\label{s:preliminaries}

We begin by introducing some notation and considering some
preliminaries on the exponential map and operator self-similarity
that are used throughout the paper.

\subsection{Some notation}

In this paper, the notation and terminology for finite-dimensional
operator theory will be preferred over their matrix analogs.
However, whenever convenient, the latter will be used.

All with respect to the field $\bbR$, $M(n)$ or $M(n,\bbR)$ is the
vector space of all $n \times n$ operators (endomorphisms),
$\mathit{GL}(n)$ or $\mathit{GL}(n,\bbR)$ is the general linear group (invertible
operators, or automorphisms), $O(n)$ is the orthogonal group of
operators $O$ such that $OO^{*} = I = O^{*}O$ (i.e., the adjoint
operator is the inverse), $\mathit{SO}(n) \subseteq O(n)$ is the special
orthogonal group of operators (rotations) with determinant equal
to 1 and $\mathit{so}(n)$ is the vector space of skew-symmetric operators
(i.e., $A^{*} = -A$).

The notation will indicate a change to the field $\bbC$. For
instance, $M(n,\bbC)$ is the vector space of complex
endomorphisms. Whenever it is said that $A \in M(n)$ has a complex
eigenvalue or eigenspace, one is considering the operator
embedding $M(n) \hookrightarrow M(n,\bbC)$. The notation
$\overline{A}$ indicates the operator whose matrix representation
is entrywise equal to the complex conjugates of those of $A$. We
will say that two endomorphisms $A, B \in M(n)$ are
\textit{conjugate} (or similar) when there exists $P \in \mathit{GL}(n)$
such that $A = P B P^{-1}$. In this case, $P$ is called a
\textit{conjugacy}. The expression
$\operatorname{diag}(\lambda_1,\ldots,\lambda_n)$ denotes the
operator whose matrix expression has the values $\lambda_1,
\ldots, \lambda_n$ on the diagonal and zeros elsewhere. The
expression $\operatorname{tr}(A)$ denotes the trace of an operator
$A \in M(n,\bbC)$. We write $f \in L^2(\bbR,M(n,\bbC))$ for a
matrix-valued function $f$ when $\operatorname{tr}\{\int_{\bbR}
f(u)^* f(u)\,\mathrm{d}u \} < +\infty$.

Throughout the paper, we set
%
\begin{equation}\label{e:D=H-(1/2)I}
D=H-(1/2)I
\end{equation}
for an operator exponent $H$. The characteristic roots of $H$ and $D$
are denoted
%
\begin{equation}\label{e:hk_dk}
h_k, d_k,
\end{equation}
respectively. Here,
%
\begin{equation}\label{e:k=1...nN}
k= 1,\ldots, N\mbox{ or }n,
\end{equation}
where $N \leq n$ is the number of different characteristic roots of
$H$.

For notational simplicity when constructing the spectral and time
domain filters, we will adopt the convention that $z^{D} = 0 \in
M(n,\bbR)$ when $z = 0$.

%

\subsection{Operator self-similar processes}
\label{ss:o.s.s._processes}

Operator self-similar (o.s.s.) processes were defined in Section \ref{s:intro}.
Any matrix $H$ for which (\ref{e:s.s.}) holds is called an
\textit{exponent} of the o.s.s. process $X$. The set of all such $H$ for
$X$ is denoted by ${\mathcal E}(X)$, which, in general, contains
more than one exponent. The non-uniqueness of the exponent~$H$
depends on the symmetry group $G_1$ of $X$, which is defined as
follows.

\begin{definition}\label{defn1}
The symmetry group of an o.s.s.\ process $X$ is the set $G_1$ of
matrices $A \in \mathit{GL}(n)$ such that
%
\begin{equation}\label{e:def_G1}
\{X(t)\}_{t \in\bbR} \stackrel{{\mathcal L}}= \{AX(t)\}_{t \in
\bbR}.
\end{equation}
\end{definition}

It turns out that the symmetry group $G_1$ is always compact,
which implies that there exists a closed subgroup ${\mathcal
O}_{0}$ of $O(n)$ such that $G_1 = W{\mathcal O}_{0}W^{-1}$, where
$W$ is a positive definite matrix (see, e.g., \cite{jurekmason1993}, Corollary 2.4.2, page 61). A process
$X$ that has maximal symmetry, that is, such that $G_1 =
WO(n)W^{-1}$, is called \textit{elliptically symmetric}.

Let $G$ be a closed (sub)group of operators. The tangent space
$T(G)$ of $G$ is the set of $A \in M(n)$ such that
\[
A = \lim_{n \rightarrow\infty} \frac{G_n - I}{d_n}
\qquad \mbox{for some } \{G_n\} \subseteq G \mbox{ and some }
0 < d_n \rightarrow0.
\]
In this sense, $T(G)$ is, in fact, a linearization of $G$ in a
neighborhood of $I$. Hudson and Mason \cite{hudsonmason1982}, Theorem 2,
shows that for any
given o.s.s.\ process $X$ with exponent $H$,
the set of exponents ${\mathcal E}(X)$ has the form ${\mathcal
E}(X) = H + T(G_1)$, where $T(G_1) = W {\mathcal L}_{0}W^{-1}$ for
the positive definite conjugacy matrix $W$ associated with $G_1$
and some subspace ${\mathcal L}_{0}$ of $\mathit{so}(n)$. Consequently, $X$
has a unique exponent if and only if $G_1$ is finite.

\section{Integral representations of OFBMs}
\label{s:integral_repres_OFBM}

Representations of OFBMs in the spectral domain are derived in
Section \ref{ss:spec_dom_repres_OFBMs}. The corresponding
representations in the time domain are given in Section
\ref{ss:timedom_repres}. The derivation of these representations
is quite different from that in the univariate case. In the latter
case, it is enough to ``guess'' the form of the spectral representation
and to verify that it gives
self-similarity and stationarity of the increments (and hence,
immediately, FBM). In the multivariate case, these representations
actually have to be derived from the properties of OFBMs, without
any guessing involved.

\subsection{Spectral domain representations}\label{ss:spec_dom_repres_OFBMs}

In Theorem \ref{t:spectral_repres_OFBM}, we establish
integral representations of OFBMs in the spectral domain.

\begin{theorem}
\label{t:spectral_repres_OFBM} Let $H \in M(n,\bbR)$ with
characteristic roots $h_k$ satisfying
%
\begin{equation}\label{e:reg_cond_char_roots_H}
0< \operatorname{Re}(h_k) < 1,\qquad  k = 1,\ldots,n.
\end{equation}
Let $\{B_{H}(t)\}_{t \in\bbR}$ be an OFBM with exponent $H$.
Then $\{B_{H}(t)\}_{t \in\bbR}$ admits the integral representation
%
\begin{equation} \label{e:spectral_repres_OFBM}
\{B_{H}(t)\}_{t \in\bbR} \stackrel{{\mathcal L}}=
\biggl\{\int_{\bbR} \frac{\mathrm{e}^{\mathrm{i}tx} - 1}{\mathrm{i}x} (x^{-D}_{+}A +
x^{-D}_{-}\overline{A}) \widetilde{B}(\mathrm{d}x) \biggr\}_{t \in\bbR}
\end{equation}
for some $A \in M(n,\bbC)$. Here, $D$ is as in (\ref{e:D=H-(1/2)I}),
%
\begin{equation}\label{e:Btilde}
\widetilde{B}(x) := \widetilde{B}_1(x) + \mathrm{i}
\widetilde{B}_{2}(x)
\end{equation}
denotes a complex-valued multivariate Brownian motion such that
$\widetilde{B}_1(-x) = \widetilde{B}_1(x)$ and
$\widetilde{B}_2(-x) = -\widetilde{B}_2(x)$, $\widetilde{B}_1$ and
$\widetilde{B}_2$ are independent and the induced random measure
$\widetilde{B}(\mathrm{d}x)$ satisfies
$E\widetilde{B}(\mathrm{d}x)\widetilde{B}(\mathrm{d}x)^{*} = \mathrm{d}x$.
\end{theorem}

\begin{pf}
For notational simplicity, set $X=B_H$. Since $X$ has stationary
increments, we have
%
\begin{equation} \label{e:X(t)-X(s)=spec_repres}
X(t) - X(s) = \int_{\bbR} \frac{\mathrm{e}^{\mathrm{i}tx} - \mathrm{e}^{\mathrm{i}sx}}{\mathrm{i}x}
\widetilde{Y}(\mathrm{d}x),
\end{equation}
where $\widetilde{Y}(\mathrm{d}x)$ is an orthogonal-increment random
measure in $\bbC^n$. The relation (\ref{e:X(t)-X(s)=spec_repres})
can be proven following the approach for the univariate case found
in  \cite{doob1953}, page\ 550, under the assumption that
$E|X(t+h)-X(t)|^2 \rightarrow0$ as $h \rightarrow0$, that is, $X$
is $L^2$-continuous at every $t$ (see also
\cite{yaglom1987}, page\ 409, and  \cite{yaglom1957},
Theorem 7). The latter assumption is satisfied in our context
because of the following. Property 2.1 in
\cite{maejimamason1994} states that, for an o.s.s.\ process $Z$
with exponent $H$, if $\inf\{\operatorname{Re}(h_k);k=1,\ldots,n\}
>0 $, then $Z(0)=0$ a.s. Thus, in view of (\ref
{e:reg_cond_char_roots_H}), $X(0)=0$ a.s.\
So, by stationarity of the increments and
continuity in law,
%
\begin{equation}\label{e:L2-cont}
X(t+h) - X(t) \stackrel{{\mathcal L}}= X(h) \stackrel{{\mathcal
L}}\rightarrow X(0) = 0,\qquad  h \rightarrow0.
\end{equation}
Therefore, by relation
(\ref{e:X(t)-X(s)=spec_repres}) and again by Property 2.1 in
\cite{maejimamason1994},
%
\begin{equation}\label{e:int(e-1/ix)Y(dx)}
X(t) = \int_{\bbR} \frac{\mathrm{e}^{\mathrm{i}tx} - 1}{\mathrm{i}x} \widetilde{Y}(\mathrm{d}x).
\end{equation}
Let $F_X(\mathrm{d}x) = E \widetilde{Y}(\mathrm{d}x) \widetilde{Y}(\mathrm{d}x)^*$ be the
multivariate spectral distribution of $\widetilde{Y}(\mathrm{d}x)$. The
remainder of the proof involves three steps:
\begin{enumerate}[(iii)]
\item[(i)] showing the existence of a spectral
density function $f_{X}(x)=F_X(\mathrm{d}x)/\mathrm{d}x$;
\item[(ii)] decorrelating the measure
$\widetilde{Y}(\mathrm{d}x)$ componentwise by finding a filter based upon
the spectral density function;
\item[(iii)] developing the form of the filter.
\end{enumerate}

Step (i). Since $X$ is o.s.s. with exponent $H$,
%
\begin{equation} \label{e:X(ct)_oss}
X(ct) \stackrel{{\mathcal L}}= c^{H} \int_{\bbR} \frac{\mathrm{e}^{\mathrm{i}tx} -
1}{\mathrm{i}x} \widetilde{Y}(\mathrm{d}x)
\end{equation}
for $c > 0$. On the other hand, through a change of variables
$x=c^{-1}v$,
%
\begin{equation} \label{e:X(ct)_change_var}
X(ct) \stackrel{{\mathcal L}}= \int_{\bbR} \frac{\mathrm{e}^{\mathrm{i}tv} - 1}{\mathrm{i}v}
c \widetilde{Y}(c^{-1}\,\mathrm{d}v).
\end{equation}
The relations (\ref{e:X(ct)_oss}) and (\ref{e:X(ct)_change_var})
provide two spectral representations for the process $\{X(ct)\}_{t
\in\bbR}$. As a consequence of the uniqueness of the spectral
distribution function of the stationary process $\{X(t) -
X(t-1)\}_{t \in\bbR}$ and of the fact that
$ |\frac{\mathrm{e}^{\mathrm{i}x}-1}{\mathrm{i}x} |^2 > 0$, $x \in\bbR\backslash\{2
\curpi k, k \in\bbZ\}$, we obtain that
\[
c^2 F_{X}(c^{-1}\,\mathrm{d}x) = c^{H}F_{X}(\mathrm{d}x)c^{H^*},\qquad  c>0.
\]
Equivalently, by a simple change of variables,
$F_X(c\,\mathrm{d}x)=c^{I-H}F_X(\mathrm{d}x)c^{(I-H)^*}$.
Thus, for $c>0$,
%
\begin{eqnarray} \label{e:F[0,c]}
\int_{(0,1]}F_{X}(c\,\mathrm{d}x) &=& F_X(0,c] = c^{I-H} F_X(0,1]
c^{(I-H)^{*}},
\\
\label{e:F[-c,0]}
\int _{(-1,0]} F_{X}(c\,\mathrm{d}x) &=& F_X(-c,0]= c^{I-H}F_X(-1,0]c^{(I-H)^*}.
\end{eqnarray}
By the explicit formula for $c^{I-H}$ in Appendix
\ref{a:Jordan_form}, each individual entry $F_{X}(0,c]_{ij}$,
$i,j=1,\ldots,n$, in the expression on the right-hand side of
(\ref{e:F[0,c]}) is either a linear combination (with complex
weights) of terms of the form
%
\begin{equation} \label{e:entries_of_F[0,c]}
\frac{(\log(c))^{l}}{l!}c^{1 - h_{q}}\frac{(\log(c))^{m}}{m!}c^{1
- \overline{h}_{k}},\qquad  q,k = 1,\ldots,n,\ l,m =
0,\ldots,n-1,
\end{equation}
or is identically zero for $c > 0$. Thus, $F_X(c)$ is
differentiable in $c$ over $(0,\infty)$ since $F_{X}(0,c]_{ij} =
F_X(c)_{ij} - F_X(0)_{ij}$. The differentiability of $F_X$ on
$(-\infty,0)$ follows from (\ref{e:F[-c,0]}) and an analogous
argument.

To finish the proof of the absolute continuity of $F_X$, it
suffices to show that $F_X$ is continuous at zero. Note that
\[
F_{X}(-c,c] = c^{I-H}F_{X}(-1,1]c^{(I-H)^*} \rightarrow0
\]
as $c \rightarrow0^+$. The limit holds because $\norm{c^{I-H}}
\rightarrow0$ as $c \rightarrow0^{+}$, where $\norm{\cdot}$ is the
matrix norm, which, in turn, follows from
\cite{maejimamason1994}, Proposition 2.1(ii),
under the assumption
that $\operatorname{Re}(h_k) < 1$, $k = 1,\ldots,n$.

Step (ii). Denote the spectral density of $X$ by $f_X$. Since
$ |\frac{1 -\mathrm{e}^{-\mathrm{i}x}}{\mathrm{i}x} |^2 f_X(x)$ is the spectral
density of the stationary process $\{X(t) - X(t-1)\}_{t \in
\bbR}$, $f_X(x)$ is a positive semidefinite Hermitian symmetric
matrix $\mathrm{d}x$-a.e. (\cite{hannan1970}, Theorem 1, page 34).
The spectral theorem yields a (unique) positive semidefinite
square root $\widehat{a}(x)$ of $f_X(x)$. Let $\widetilde{B}(x)$
be a complex-valued multivariate Brownian motion, as in the
statement of the theorem. $X$ can then also be represented as
%
\begin{equation}\label{e:int(e-1/ix)a(x)B(dx)}
X(t) \stackrel{{\mathcal L}}= \int_{\bbR} \frac{\mathrm{e}^{\mathrm{i}tx} - 1}{\mathrm{i}x}
\widehat{a}(x) \widetilde{B}(\mathrm{d}x)
\end{equation}
because
\[
E (\widehat{a}(x) \widetilde{B}(\mathrm{d}x) \widetilde{B}(\mathrm{d}x)^{*}
\widehat{a}(x)^{*}) =\widehat{a}(x)^2 \,\mathrm{d}x = f_X(x) \,\mathrm{d}x = F_X(\mathrm{d}x)
\]
and the processes on both sides of
(\ref{e:int(e-1/ix)a(x)B(dx)}) are
Gaussian and real-valued.

Step (iii). By using operator self-similarity and arguing as in
step (i), the relation (\ref{e:int(e-1/ix)a(x)B(dx)}) implies
that, for every $c
> 0$,
%
\begin{equation}\label{e:a(x)a(x)*_scaling}
\widehat{a}(x)\widehat{a}(x)^{*}=c^{-D}\widehat{a} \biggl(\frac{x}{c}
\biggr)\widehat{a} \biggl(\frac{x}{c} \biggr)^*c^{-D^*}
\qquad \mbox{$\mathrm{d}x$-a.e.}
\end{equation}
By Fubini's theorem, the relation (\ref{e:a(x)a(x)*_scaling}) also
holds $\mathrm{d}x\,\mathrm{d}c$-a.e.

Consider $x>0$. A change of variables leads to
\[
\widehat{a}(x)\widehat{a}(x)^{*}=x^{-D}v^{D}\widehat{a}(v)\widehat
{a}(v)^*v^{D^*}x^{-D^*}\qquad
\mbox{$\mathrm{d}x \,\mathrm{d}v$-a.e.}
\]
Thus, one can choose $v_+ > 0$ such that
%
\begin{equation}\label{e:squares_filters}
\widehat{a}(x)\widehat{a}(x)^{*}=x^{-D}v^{D}_+
\widehat{a}(v_+)\widehat{a}(v_+)^*v^{D^*}_+x^{-D^*}
\qquad \mbox{$\mathrm{d}x$-a.e.}
\end{equation}
This means, in particular, that if we set
\[
\widehat{\alpha}_+(x) = x^{-D}v^{D}_+ \widehat{a}(v_+)
\]
for $\mathrm{d}x$-a.e.\ $x > 0$, then
$\widehat{\alpha}_+(x)\widehat{\alpha}_+(x)^*=f_X(x)$ on the same
domain.

By again considering the stationary process $\{X(t) - X(t-1)\}_{t
\in\bbR}$ and applying \cite{hannan1970}, Theorem~3, page 41, one can show that $f_X$ is a Hermitian
function. Thus,
\[
\widehat{a}(-x)\widehat{a}(-x)^{*}=f_{X}(-x) = \overline{f_{X}(x)}
= x^{-D}v^{D}_+
\overline{\widehat{a}(v_+)\widehat{a}(v_+)^*} v^{D^*}_+x^{-D^*}
\qquad \mbox{$\mathrm{d}x$-a.e.}
\]
Hence, for $x < 0$, we can set
\[
\widehat{\alpha}_{-}(x)=(-x)^{-D}_+ v^{D}_+
\overline{\widehat{a}(v_+)}
\]
and, for $x \in\bbR$, we have
\[
\widehat{\alpha}(x)=x^{-D}_+ v^{D}_+
\widehat{a}(v_+)+x^{-D}_{-}v^{D}_{+}
\overline{\widehat{a}(v_+)}
\qquad \mbox{$\mathrm{d}x$-a.e.},
\]
where $\widehat{\alpha}(x)\widehat{\alpha}(x)^*=f_X(x)$ $\mathrm{d}x$-a.e.\
Therefore, we can use $\widehat{\alpha}$ in place of $\widehat{a}$
in the spectral representation of $X$, which establishes relation
(\ref{e:spectral_repres_OFBM}).
\end{pf}

\begin{remark}\label{r:example_proper_with_deficient_rank_A} The
invertibility of $A$ in relation
(\ref{e:spectral_repres_OFBM}) is not a requirement for the
process to be proper (compare with
\cite{bahadoranbenassidebicki2003}, page 9). In the Gaussian case, properness
is equivalent to $EX(t)X(t)^*$ being a full rank matrix for all $t
\neq0$.

A simple example would be that of a bivariate OFBM whose spectral
representation has matrix parameters $D = d I$, $0 < d
< 1/2$, and $A$ set to the (unique) non-negative square root of
\[
A^2 = \pmatrix{ 1 & i \cr \overline{i} & 1},
\]
which is rank-deficient. Let
\[
g(t) = \int^{\infty}_{0} \biggl|\frac{\mathrm{e}^{\mathrm{i}tx}-1}{\mathrm{i}x} \biggr|^2
|x|^{-2d}\,\mathrm{d}x,
\]
which is strictly positive for all $t \neq0$. In this case,
\begin{eqnarray*}
EX(t)X(t)^* &=& \int_{\bbR} \biggl|\frac{\mathrm{e}^{\mathrm{i}tx}-1}{\mathrm{i}x} \biggr|^2
|x|^{-2d}\bigl(A^2 1_{\{x \geq0\}}+\overline{A^2} 1_{\{x < 0\}}\bigr)\,\mathrm{d}x
\\
&=&g(t) \pmatrix{ 2 & i+\overline{i} \cr i+\overline{i}
& 2} = 2g(t)I.
\end{eqnarray*}
\end{remark}

Theorem \ref{t:spectral_repres_OFBM} shows that an OFBM is
characterized by a (potentially non-unique) o.s.s.\ exponent $H$
and a matrix $A$. For the sake of simplicity, we will continue to
use the notation $B_H$ instead of the (more correct) notation
$B_{H,A}$.

\begin{remark}\label{rem:no_OFBM_with_H=Jordan_h=1}
As a consequence of Maejima and Mason
\cite{maejimamason1994}, Corollary 2.1, the characteristic roots
$h_k$ of the
exponent $H$ of an OFBM must satisfy $\operatorname{Re}(h_k)\leq1$,
$k=1,\ldots,n$. However, the extension of the definition of  OFBMs to
the case of $H$ with at least one characteristic root $h_k$
satisfying $\operatorname{Re}(h_k)= 1$ can be subtle. In Proposition
\ref{p:no_OFBM_with_H=Jordan_h=1}, it is shown that there does not
exist an OFBM with exponent
\[
H = \pmatrix{ 1 & 0\cr
1 & 1
}
\]
whose characteristic roots are $h_1 = h_2 = 1$. (More precisely,
it is shown that a Gaussian, $H$-o.s.s.\ process $X =
(X_1,X_2)^{*}$ with stationary increments is necessarily such that
$X_1(t)=0$ and $X_2(t)=tY$ a.s.\ for a Gaussian variable $Y$ and
hence that it cannot be proper.)
\end{remark}

\subsection{Time domain representations}\label{ss:timedom_repres}

Our next goal is to provide integral representations of OFBMs in
the time domain, which is done in Theorem \ref{t:time domain_OFBM}.
The key technical step in the proof is the calculation of
the (entrywise) Fourier transform of the kernels
%
\begin{equation}\label{e:OFBM_series_time_domain_kernel}
(t-u)^{D}_{\pm}-(-u)^{D}_{\pm} = \exp\bigl(\log(t-u)_{\pm}D\bigr)-
\exp\bigl(\log(-u)_{\pm}D\bigr),
\end{equation}
which are the multivariate analogs of the corresponding
univariate FBM time domain kernels. It is natural and convenient
to carry out this step in the framework of the so-called primary
matrix functions. The latter allows one to naturally define matrix
analogs $f(D)$, $D \in M(n,\bbR)$, of univariate functions
$f(d)$, $d \in\bbR$, and to say when two such matrix-valued
functions are equal based on their univariate counterparts.

For the reader's convenience, we recall the definition of primary
matrix functions (more details and properties can be found in
\cite{hornjohnson1991}, Sections 6.1 and 6.2). Let
$\Lambda\in M(n,\bbC)$ with minimal polynomial
%
\begin{equation}\label{e:q(t)}
q_{\Lambda}(z) = (z - \lambda_1)^{r_1}\cdots(z - \lambda_N)^{r_N},
\end{equation}
where $\lambda_1,\ldots,\lambda_N$ are pairwise distinct and
$r_k \geq1$ for $k=1,\ldots,N$, $N \leq n$. We denote by $\Lambda
=PJP^{-1}$ the Jordan decomposition of
$\Lambda$, where $J$ is in Jordan canonical form with the Jordan
blocks $J_{\lambda_1},\ldots, J_{\lambda_N}$ on the diagonal.

Let $U \subseteq\bbC$ be an open set. Given a function $h\dvtx U
\rightarrow\bbC$ and some $\Lambda\in M(n,\bbC)$ as above,
consider the following conditions: (M1) $\lambda_k \in U$, $k=1,\ldots
,N$; (M2) if
$r_k > 1$, then $h(z)$ is analytic in a neighborhood $U_k \ni
\lambda_k$, where $U_k \subseteq U$.
Let ${\mathcal M}_h =\{ \Lambda\in M(n,\bbC)$; conditions
(M1) and (M2) hold at the characteristic roots
$\lambda_1,\ldots,\lambda_N$ of $\Lambda$\}.
We now define the primary matrix function $h(\Lambda)$ associated
with the scalar-valued stem function $h(z)$.

\begin{definition}\label{d:matrix_function}
The primary matrix function $h:{\mathcal M}_h \rightarrow
M(n,\bbC)$ is defined as
\[
h(\Lambda) = Ph(J)P^{-1}=P \pmatrix{
h(J_{\lambda_1})
& \ldots& 0\cr
\vdots& \ddots& \vdots\cr
0 & \ldots& h(J_{\lambda_N})
}P^{-1},
\]
where
\[
h(J_{\lambda_k})= \pmatrix{ h(\lambda_k) & 0 & \ldots& 0 \cr
h'(\lambda_k) & h(\lambda_k) & \ddots& 0\cr
\vdots& \ddots& \ddots& \vdots\vspace*{3pt}\cr
\displaystyle\frac{h^{(r_k-1)}(\lambda_k)}{(r_k -1)!} & \ldots & h'(\lambda_k) &
h(\lambda_k) }.
\]
\end{definition}

The following technical result is proved in Appendix \ref{ss:FT}. The functions
$(t-u)^{D}_{\pm}$, $\Gamma(D+I)$, $|x|^{-D}$, $\mathrm{e}^{\mp\operatorname{sign}(x)\mathrm{i}
\curpi D/2}$ appearing in the result below are all primary matrix functions.
The same interpretation
is also adopted throughout the rest of the paper, for example, with
functions $\sin(\curpi D/2)$,
$\cos(\curpi D/2)$ appearing in Theorem \ref{t:time domain_OFBM}.
(It should also be noted, in particular,
that the definition of the matrix exponential based on a series is
equivalent to that based
on primary matrix functions.)

\begin{proposition}\label{p:FT_timedom_filter}
Under (\ref{e:reg_cond_char_roots_H}) and condition
(\ref{e:ch_roots_neq_1/2}) in Theorem \ref{t:time domain_OFBM},
%
\begin{equation}\label{e:FT_timdom_filter}
\int_{\bbR} \mathrm{e}^{\mathrm{i}ux} \bigl((t - u)^{D}_{\pm} - (- u)^{D}_{\pm}
\bigr)\,\mathrm{d}u
 = \frac{\mathrm{e}^{\mathrm{i}tx}-1}{\mathrm{i}x}|x|^{-D}\Gamma(D+I)\mathrm{e}^{\mp
\operatorname{sign}(x)\mathrm{i} \curpi D/2}.
\end{equation}
\end{proposition}

Next, we construct time domain representations for OFBMs, which provides
the main result in this section. Further comments about the result
can be found after the proof.

\begin{theorem}\label{t:time domain_OFBM}
Let $\{B_{H}(t)\}_{t \in\bbR}$ be an OFBM with o.s.s.\ exponent
$H$ having the spectral representation
(\ref{e:spectral_repres_OFBM}) with $A = A_1 + \mathrm{i} A_2$, where $A_1,
A_2 \in M(n,\bbR)$.
\begin{enumerate}[(ii)]
\item[(i)] Suppose that $H \in M(n,\bbR)$ has characteristic
roots satisfying (\ref{e:reg_cond_char_roots_H}) and
%
\begin{equation}\label{e:ch_roots_neq_1/2}
\operatorname{Re}(h_k) \neq\frac{1}{2},\qquad  k = 1,\ldots,n.
\end{equation}
There are then $M_{+}, M_{-} \in M(n,\bbR)$ such that
\begin{eqnarray}\label{e:time_repres_OFBM}
&&\{B_{H}(t)\}_{t \in\bbR} \nonumber\\[-8pt]\\[-8pt]
&&\quad \stackrel{{\mathcal L}}= \biggl\{\int_{\bbR} \bigl(\bigl((t -
u)^{D}_{+} - (- u)^{D}_{+}\bigr)M_{+} + \bigl((t - u)^{D}_{-} - (- u)^{D}_{-}\bigr)M_{-}
\bigr) B(\mathrm{d}u) \biggr\}_{t \in\bbR},\nonumber\qquad
\end{eqnarray}
where $\{B(u)\}_{u \in\bbR}$ is a vector-valued process
consisting of independent Brownian motions and such that
$EB(\mathrm{d}u)B(\mathrm{d}u)^*=\mathrm{d}u$. Moreover, the matrices $M_{+}$, $M_{-}$ can be
taken as
%
\begin{equation}\label{e:M_pm}
M_{\pm} = \sqrt{\frac{\curpi}{2}} \biggl(\sin\biggl(\frac{\curpi
D}{2} \biggr)^{-1}\Gamma(D+I)^{-1}A_1 \pm\cos\biggl(\frac{\curpi
D}{2} \biggr)^{-1}\Gamma(D+I)^{-1}A_2 \biggr).
\end{equation}

\item[(ii)] Suppose that $H = (1/2)I$. There then exist $M,N \in
M(n,\bbR)$ such that
%
\begin{equation}\label{e:time domain_OFBM_D=0}
\{B_H(t)\}_{t \in\bbR} \stackrel{{\mathcal L}}= \biggl\{\int_{\bbR}
\biggl(\bigl(\operatorname{sign}(t-u)-\operatorname{sign}(-u)\bigr)M+
\log\biggl(\frac{|t-u|}{|u|} \biggr)N \biggr) B(\mathrm{d}u) \biggr\}_{t \in\bbR},
\end{equation}
where $\{B(u)\}_{u \in\bbR}$ is as in (\ref{e:time_repres_OFBM}).
Moreover, the matrices $M$, $N$ can be taken as
%
\begin{equation}\label{e:M_N}
M = \sqrt{\frac{\curpi}{2}} A_1,\qquad  N =
-\sqrt{\frac{2}{\curpi}} A_2.
\end{equation}
\end{enumerate}
\end{theorem}

\begin{pf}
(i) Denote the process on the right-hand side of
(\ref{e:time_repres_OFBM}) by $X_{H}$. By using the Jordan
decomposition of $D$, it is easy to show that $X_H$ is
well defined. It suffices to show that there are $M_{\pm}$ such
that the covariance structure of $X_H$ matches that of the OFBM
$B_H$ given by its spectral representation
(\ref{e:spectral_repres_OFBM}) with $A = A_1 + \mathrm{i} A_2$. By using
the Plancherel identity, note first that\looseness=1
\begin{eqnarray*}
&&E X_{H}(s)X_{H}(t)^*\\
&&\quad = \frac{1}{2\curpi}\int_{\bbR}
\frac{(\mathrm{e}^{\mathrm{i}sx}-1)(\mathrm{e}^{-\mathrm{i}tx}-1)}{|x|^2} \bigl(|x|^{-D}
\Gamma(D+I)\bigl(\mathrm{e}^{-\operatorname{sign}(x)\mathrm{i} \curpi D/2}M_{+} + \mathrm{e}^{
\operatorname{sign}(x)\mathrm{i} \curpi D/2}M_{-}\bigr) \bigr)
\\
&&\hphantom{\frac{1}{2\curpi}\int_{\bbR}}\qquad {}\times \bigl(\bigl(M^{*}_{+} \mathrm{e}^{ \operatorname{sign}(x)\mathrm{i} \curpi D^*/2}
+M^{*}_{-} \mathrm{e}^{-\operatorname{sign}(x)\mathrm{i} \curpi
D^{*}/2}\bigr)\Gamma(D+I)^{*}|x|^{-D^{*}} \bigr)\,\mathrm{d}x.
\end{eqnarray*}
Meanwhile, for $B_H$, we have
%
\begin{equation} \label{e:cov_struct_BH(t)}
E B_H(s)B_H(t)^* =
\int_{\bbR}\frac{(\mathrm{e}^{\mathrm{i}sx}-1)(\mathrm{e}^{-\mathrm{i}tx}-1)}{|x|^2}
(x^{-D}_{+}AA^{*}x^{-D^*}_{+}+x^{-D}_{-}\overline{AA^{*}}x^{-D^*}_{-})\,\mathrm{d}x.
\end{equation}
Thus, by using the relation $\mathrm{e}^{\mathrm{i} \Theta} = \cos(\Theta) + \mathrm{i}
\sin(\Theta)$, $\Theta\in M(n)$, it is sufficient to find
$M_{\pm} \in M(n,\bbR)$ such that
\begin{eqnarray} \label{e:AA*_on_time_domain}
AA^{*} &=& \frac{1}{2 \curpi}\Gamma(D+I)(\mathrm{e}^{-\mathrm{i}\curpi D/2}M_{+} + \mathrm{e}^{\mathrm{i}\curpi D/2}M_{-}) \nonumber
\\
&&{}\times(M^{*}_{+}\mathrm{e}^{\mathrm{i}\curpi D^*/2} + M^{*}_{-}\mathrm{e}^{-\mathrm{i}\curpi D^{*}/2})\Gamma(D+I)^{*}\nonumber
\\
&=& \frac{1}{2 \curpi} \Gamma(D+I) \biggl( \sin\biggl(\frac{\curpi D}{2} \biggr)(M_{+}-M_{-})(M^*_{+}-M^*_{-})\sin\biggl(\frac{\curpi
D^*}{2} \biggr)\nonumber\\[-8pt]\\[-8pt]
&&\hphantom{\frac{1}{2 \curpi} \Gamma(D+I) \biggl( }{}+ \cos\biggl(\frac{\curpi
D}{2} \biggr)(M_{+}+M_{-})(M^*_{+}+M^*_{-})\cos\biggl(\frac{\curpi
D^*}{2} \biggr)\nonumber
\\
&&\hphantom{\frac{1}{2 \curpi} \Gamma(D+I) \biggl( }{}+ \mathrm{i} \biggl(\cos\biggl(\frac{\curpi
D}{2} \biggr)(M_{+}+M_{-})(M^*_{+}-M^*_{-})\sin\biggl(\frac{\curpi
D^*}{2} \biggr)\nonumber
\\
&&\hphantom{\frac{1}{2 \curpi} \Gamma(D+I) \biggl( {}+ \mathrm{i} \biggl(}{}-\sin\biggl(\frac{\curpi
D}{2} \biggr)(M_{+}-M_{-})(M^*_{+}+M^*_{-})\cos\biggl(\frac{\curpi
D^*}{2} \biggr) \biggr) \biggr)\Gamma(D+I)^{*}.\nonumber
\end{eqnarray}
On the other hand,
%
\begin{equation}\label{e:AA*}
AA^* = (A_1 A^{*}_1 + A_2 A^{*}_2) + \mathrm{i}(A_2 A^{*}_1 - A_1 A^{*}_2).
\end{equation}
By comparing (\ref{e:AA*}) and (\ref{e:AA*_on_time_domain}), a
natural way to proceed is to consider $M_+$ and $M_-$ as solutions
of the system
\begin{eqnarray}\label{e:M+_and_M-_from_A1_and_A2}
A_1 &=& \frac{1}{\sqrt{2
\curpi}} \Gamma(D+I)\sin\biggl(\frac{\curpi
D}{2} \biggr)(M_{+}-M_{-}),\nonumber\\[-8pt]\\[-8pt]
A_2 &=& \frac{1}{\sqrt{2
\curpi}} \Gamma(D+I)\cos\biggl(\frac{\curpi
D}{2} \biggr)(M_{+}+M_{-}).\nonumber
\end{eqnarray}
By assumption (\ref{e:ch_roots_neq_1/2}), $\sin(\frac{\curpi
D}{2} )$, $\cos(\frac{\curpi D}{2} )$ and $\Gamma(D+I)$ are
invertible,
and we obtain the solution given by (\ref{e:M_pm}).

(ii) In this case, one can readily compute the inverse Fourier
transform of the integrand in (\ref{e:spectral_repres_OFBM}), that
is (up to $(2\curpi)^{-1}$),
\begin{eqnarray*}
\hspace*{-4.5pt}&&\int_{\bbR}\mathrm{e}^{-\mathrm{i}ux} \biggl(\frac{\mathrm{e}^{\mathrm{i}tx}-1}{\mathrm{i}x} \biggr) \bigl(1_{\{x > 0
\}}A + 1_{\{x < 0 \}}\overline{A} \bigr)\,\mathrm{d}x
\\
\hspace*{-4.5pt}&&\quad =\int_{\bbR} \biggl(\frac{\cos((t-u)x)-\cos((t-u)x)+\mathrm{i}(\sin((t-u)x)+\sin
(ux))}{\mathrm{i}x} \biggr)
\bigl(1_{\{x > 0 \}} A + 1_{\{x < 0 \}}\overline{A}\bigr)\,\mathrm{d}x.
\end{eqnarray*}
As shown in Appendix \ref{a:useful_integ}, this becomes
\[
-2\log\biggl(\frac{|t-u|}{|u|} \biggr)A_2 +
\bigl(\operatorname{sign}(t-u)-\operatorname{sign}(-u)\bigr)\curpi A_1 .
\]
Then, by considering second moments and using Plancherel's
identity, representation (\ref{e:time domain_OFBM_D=0}) holds with
$M =(2 \curpi)^{-1/2} \curpi A_1 $ and $N= (2 \curpi)^{-1/2} (-2) A_2$. It
is well defined because the integrand comes from the inverse
Fourier transform of a square-integrable function and hence is
also square-integrable.
\end{pf}

\begin{remark}
Note that the invertibility of $M$ or $N$ in
(\ref{e:time_repres_OFBM}) is not a requirement for the process to
be proper. A simple example would be that of a bivariate OFBM
$B_H$ whose time domain representation (\ref{e:time_repres_OFBM})
has matrix parameters $H = h I$, $h \in
(0,1)\backslash\{1/2\}$,
\[
M = \pmatrix{ 1 & 0\cr
0 & 0},
\qquad
N = \pmatrix{ 0 & 0\cr
0 & 1}.
\]
The two components of $B_H$ are two independent (univariate) FBMs
with exponent $h$. Thus, $B_H$ is proper.
\end{remark}

\begin{example}\label{ex:time domain_D=Jordan_and_eig=0}
When (\ref{e:ch_roots_neq_1/2}) does not hold and $H \neq(1/2)I$, the general
form of time domain representations can be quite intricate. For
example, with
\[
D = \pmatrix{
0 & 0\cr
1 & 0
}\qquad
\left(H = \pmatrix{
1/2 & 0\cr
1 & 1/2
} \right),
\]
the calculation of the inverse Fourier transform (up to $(2
\curpi)^{-1}$)
%
\begin{equation}\label{e:IFT_D=Jordan_and_eig=0}
\int_{\bbR}\mathrm{e}^{-\mathrm{i}ux} \biggl(\frac{\mathrm{e}^{\mathrm{i}tx}-1}{\mathrm{i}x}
\biggr)(x^{-D}_{+}A+x^{-D}_{-}\overline{A})\,\mathrm{d}x
\end{equation}
in Appendix \ref{a:useful_integ} shows that $B_H$ has the time domain
representation
%
\begin{equation}\label{e:time domain_D=Jordan_and_eig=0}
\{B_H(t)\}_{t \in\bbR} \stackrel{{\mathcal L}}=
\biggl\{\int_{\bbR} \bigl(f_1(t,u)M+f_2(t,u)N \bigr)B(\mathrm{d}u) \biggr\}_{t \in
\bbR},
\end{equation}
where $M = \sqrt{\frac{\curpi}{2}} A_1$,
$N= -\sqrt{\frac{2}{\curpi}} A_2$,
\begin{eqnarray*}
f_1(t,u) &=& \pmatrix{
\operatorname{sign}(t-u)-\operatorname{sign}(-u) & 0\cr
(C+\log|t-u|)\operatorname{sign}(t-u)-(C+\log|u|)\operatorname{sign}(-u)
& \operatorname{sign}(t-u)-\operatorname{sign}(-u)
},
\\
f_2(t,u) &=& \pmatrix{
\displaystyle\log\biggl(\frac{|t-u|}{|u|} \biggr) & 0\cr
\displaystyle\log\biggl(\frac{|t-u|}{|u|} \biggr) \biggl(C+\frac{1}{2}\log(|t-u||u|) \biggr)
& \displaystyle\log\biggl(\frac{|t-u|}{|u|} \biggr)
},
\end{eqnarray*}
where $C$ is Euler's constant. Note that, without taking the
Fourier transform of (\ref{e:time domain_D=Jordan_and_eig=0}), it
is by no means obvious why its right-hand side has stationary
increments and is o.s.s.
\end{example}

\section{Conditions for properness}\label{s:properness}

We now provide sufficient conditions for a process with spectral
and time domain representations~(\ref{e:spectral_repres_OFBM}) and
(\ref{e:time_repres_OFBM}), respectively, to be proper and, thus,
to be an OFBM.
\begin{proposition}\label{p:conditions_properness}
Let $\{X(t)\}_{t \in\bbR}$ be a process with spectral
domain representation (\ref{e:spectral_repres_OFBM}), where the
characteristic roots of $H$ satisfy
(\ref{e:reg_cond_char_roots_H}). If $\operatorname{Re}(AA^*)$ is a
full rank matrix, then $\{X(t)\}_{t \in\bbR}$ is proper (i.e., it
is an OFBM).
\end{proposition}
\begin{pf}
We must show that
\[
EX(t)X(t)^* = \int_{\bbR} \biggl|
\frac{\mathrm{e}^{\mathrm{i}tx}-1}{\mathrm{i}x} \biggr|^2(x^{-D}_{+}AA^*x^{-D^*}_{+}+x^{-D}_{-}\overline
{AA^*}x^{-D^*}_{-})\,\mathrm{d}x,
\qquad t \neq0,
\]
is a full rank matrix. For simplicity, let $\mathrm{d}\mu(x)= |
\frac{\mathrm{e}^{\mathrm{i}tx}-1}{\mathrm{i}x} |^2\,\mathrm{d}x$. Then
\begin{eqnarray*}
EX(t)X(t)^* &=& \int_{\bbR} x^{-D}_{+}AA^*x^{-D^*}_{+} \,\mathrm{d}\mu(x)
+\int_{\bbR}x^{-D}_{+}\overline{AA^*}x^{-D^*}_{+} \,\mathrm{d}\mu(x)
\\
&=& 2\int_{\bbR} x^{-D}_{+}\operatorname{Re}(AA^*)x^{-D^*}_{+}\,\mathrm{d}\mu(x).
\end{eqnarray*}
The matrix $\int_{\bbR}
x^{-D}_{+}\operatorname{Re}(AA^*)x^{-D^*}_{+} \,\mathrm{d}\mu(x)$ is Hermitian
positive semidefinite. Moreover, for any $v \in\bbC^n \backslash
\{0\}$,
\[
v^* \biggl( \int_{\bbR} x^{-D}_{+}\operatorname{Re}(AA^*)x^{-D^*}_{+}
\,\mathrm{d}\mu(x) \biggr) v > 0,
\]
where the strict inequality follows from the fact that $(v^*
x^{-D}_{+})\operatorname{Re}(AA^*)(x^{-D^*}_{+} v) > 0$ for all \mbox{$x >
0$}, the latter being a consequence of the invertibility of
$x^{-D}_{+}$ and the assumption that
$\operatorname{Re}(AA^*)$ has full rank.
\end{pf}

Based on Proposition \ref{p:conditions_properness}, we can easily
obtain conditions for properness based on time domain parameters.
Consider a process $\{X(t)\}_{t \in\bbR}$ with time domain
representation (\ref{e:time_repres_OFBM}), where the
characteristic roots of $H$ satisfy
(\ref{e:reg_cond_char_roots_H}) and (\ref{e:ch_roots_neq_1/2}). If
\[
M_+ + M_-,\qquad  M_+ - M_-
\]
are full rank matrices, then $\{X(t)\}_{t \in\bbR}$ is proper
(i.e., it is an OFBM).

\begin{remark}
$\operatorname{Re}(AA^*)$ having full rank does not imply that
$AA^*$ has full rank since $\mathrm{i}(A_2 A^{*}_{1}-A_1A^{*}_2)$ may have
negative eigenvalues. Also, note that $\operatorname{Re}(AA^*)$ being a full
rank matrix is not a
necessary condition for properness. For example, consider the
process $\{X(t)\}_{t \in\bbR}$ with spectral representation
(\ref{e:spectral_repres_OFBM}), where
\[
AA^* = \pmatrix{
1 & 2\cr
2 & 4
},\qquad
H = \pmatrix{
h_1 & 0\cr
0 & h_2
},\qquad
h_1,h_2 \in(0,1).
\]
Then
\begin{eqnarray*}
EX(t)X(t)^* &=& \int_{\bbR} \biggl|\frac{\mathrm{e}^{\mathrm{i}tx}-1}{\mathrm{i}x} \biggr|^2
\pmatrix{
|x|^{-2(h_1-1/2)} & 2|x|^{-((h_1-1/2)+(h_2-1/2))}\cr
2|x|^{-((h_1-1/2)+(h_2-1/2))} & 4|x|^{-2(h_2-1/2)}
}\,\mathrm{d}x
\\
&=& \pmatrix{
|t|^{2 h_1} C_{2}(h_1)^2 &\displaystyle 2 |t|^{h_1 + h_2} C_{2}\biggl(\frac
{h_1+h_2}{2}\biggr)^2\cr
\displaystyle 2 |t|^{h_1 + h_2} C_{2}\biggl(\frac{h_1+h_2}{2}\biggr)^2 & 4 |t|^{2 h_2}
C_{2}(h_2)^{2}
},
\end{eqnarray*}
where
%
\begin{equation}\label{e:C2(H)}
C_{2}(h)^2 = \frac{\curpi}{h \Gamma(2h) \sin(h\curpi)}
\end{equation}
(see, e.g., \cite{samorodnitskytaqqu1994}, page\ 328). Therefore,
$\operatorname{det}(EX(t)X(t)^*)=0$ if and only if
\[
C_{2}(h_1)^2C_{2}(h_2)^2 =
\biggl(C_{2} \biggl(\frac{h_1+h_2}{2} \biggr)^2 \biggr)^2,
\]
which generally does not hold.

\end{remark}

\section{Time-reversibility of OFBMs}\label{s:time_reversibility}

We shall provide here conditions for an OFBM to be time-reversible.
Recall that
a process $X$ is said to be time-reversible if
%
\begin{equation}\label{e:isotropic_processes}
\{X(t)\}_{t \in\bbR} \stackrel{{\mathcal L}}= \{X(-t)\}_{t \in
\bbR}.
\end{equation}
When $X$ is a zero-mean multivariate Gaussian stationary process,
(\ref{e:isotropic_processes}) is equivalent to
\[
EX(s)X(t)^* = EX(-s)X(-t)^*,\qquad  s,t \in\bbR,
\]
which, in turn, is equivalent to
\[
EX(s)X(t)^* = EX(t)X(s)^*,\qquad  s,t \in\bbR.
\]
The next proposition provides necessary and sufficient conditions
for time-reversibility in the case of Gaussian processes with
stationary increments. It is stated without proof since the
latter is elementary.

\begin{proposition}\label{p:time-revers_s.i._processes}
Let $X$ be a Gaussian process with stationary increments and
spectral representation
\[
\{X(t)\}_{t \in\bbR} \stackrel{{\mathcal L}}=
\biggl\{\int_{\bbR} \frac{\mathrm{e}^{\mathrm{i}tx}-1}{\mathrm{i}x}
\widetilde{Y}(\mathrm{d}x) \biggr\}_{t \in\bbR},
\]
where $\widetilde{Y}(\mathrm{d}x)$ is an orthogonal-increment random
measure in $\bbC^n$. The following statements are equivalent:
\begin{longlist}[(iii)]
\item[(i)] $X$ is time-reversible;
\item[(ii)] $E
\widetilde{Y}(\mathrm{d}x)\widetilde{Y}(\mathrm{d}x)^* = E
\widetilde{Y}(-\mathrm{d}x)\widetilde{Y}(-\mathrm{d}x)^*$;
\item[(iii)]
$EX(s)X(t)^* = EX(t)X(s)^*$, $s,t \in\bbR$.
\end{longlist}
\end{proposition}

The following result on time-reversibility of OFBMs is a direct
consequence of Proposition~\ref{p:time-revers_s.i._processes}.

\begin{theorem}\label{t:time-revers_OFBM}
Let $\{B_H(t)\}_{t \in\bbR}$ be an OFBM with exponent $H$ and spectral
representation
(\ref{e:spectral_repres_OFBM}). Let $A = A_1 + \mathrm{i} A_2$, where $A_1, A_2
\in M(n,\bbR)$. Then $B_H$ is time-reversible if and
only if
%
\begin{equation}\label{e:A2A^*1=A1A^*2}
AA^* = \overline{AA^*} \quad \mbox{or}\quad  A_2 A^*_1 = A_1 A^*_2.
\end{equation}
\end{theorem}
\begin{pf}
From Proposition \ref{p:time-revers_s.i._processes}(ii),
time-reversibility is equivalent to
\begin{eqnarray*}
&&E \bigl((x^{-D}_{+}A+x^{-D}_{-}\overline{A})\widetilde{B}(\mathrm{d}x)\widetilde{B}(\mathrm{d}x)^*
(A^*x^{-D^*}_{+} + \overline{A^*}x^{-D^*}_{-}) \bigr)
\\
&&\quad =E \bigl((x^{-D}_{-}A+x^{-D}_{+}\overline{A})\widetilde{B}(-\mathrm{d}x)\widetilde{B}(-\mathrm{d}x)^*
(A^*x^{-D^*}_{-} + \overline{A^*}x^{-D^*}_{+}) \bigr)
\end{eqnarray*}
or
\[
x^{-D}_{+}AA^*x^{-D}_{+}+x^{-D}_{-}\overline{AA^*}x^{-D^*}_{-} =
x^{-D}_{-}AA^*x^{-D}_{-}+x^{-D}_{+}\overline{AA^*}x^{-D^*}_{+}\qquad
\mbox{$\mathrm{d}x$-a.e.}
\]
Since $|x|^{D}$ is invertible for $x > 0$, this is equivalent to
(\ref{e:A2A^*1=A1A^*2}).
\end{pf}
\begin{corollary}\label{c:time_revers_timedom_condns}
Let $\{B_H(t)\}_{t \in\bbR}$ be an OFBM with time domain
representation given by
(\ref{e:time_repres_OFBM}) and exponent $H$ satisfying
(\ref{e:reg_cond_char_roots_H}) and (\ref{e:ch_roots_neq_1/2}).
Then $B_H$ is time-reversible if and only if
\begin{eqnarray} \label{e:time_revers_timedom_condns}
&&\cos\biggl(\frac{\curpi D}{2} \biggr)(M_{+}+M_{-})(M^{*}_{+}-M^{*}_{-}) \sin\biggl(\frac
{\curpi
D^*}{2} \biggr)
\nonumber\\[-8pt]\\[-8pt]
&&\quad =\sin\biggl(\frac{\curpi
D}{2} \biggr)(M_{+}-M_{-}) (M^{*}_{+}+M^{*}_{-}) \cos\biggl(\frac{\curpi D^*}{2} \biggr).\nonumber
\end{eqnarray}
\end{corollary}
\begin{pf}
As in the proof of Theorem \ref{t:time domain_OFBM}, under
(\ref{e:ch_roots_neq_1/2}), the matrices $\sin(\curpi D/2)$, $\cos(\curpi
D/2)$ and $\Gamma(D+I)$ are invertible and, thus, by using
(\ref{e:M+_and_M-_from_A1_and_A2}), one can equivalently
re-express condition~(\ref{e:A2A^*1=A1A^*2}) as (\ref
{e:time_revers_timedom_condns}).
\end{pf}

A consequence of Theorem \ref{t:time-revers_OFBM} is that
time-irreversible OFBMs can only emerge in the multivariate context
since condition (\ref{e:A2A^*1=A1A^*2})
is always satisfied in the univariate context. Another elementary
consequence of Proposition
\ref{p:time-revers_s.i._processes} is the following result, which
partially justifies the interest in time-reversibility in the case
of OFBMs.

\begin{proposition}
Let $\{B_{H}(t)\}_{t \in\bbR}$ be an OFBM with $H$ satisfying
(\ref{e:reg_cond_char_roots_H}). If $\{B_{H}(t)\}_{t \in\bbR}$ is
time-reversible, then its covariance structure is given by the
function
%
\begin{eqnarray}\label{e:cov_func_timerevers_OFBM}
EB_H(s)B_H(t)^* &=&
\frac{1}{2} \bigl(|t|^{H}\Gamma(1,1)|t|^{H^*}+|s|^{H}\Gamma
(1,1)|s|^{H^*}\nonumber\\[-8pt]\\[-8pt]
&&\hphantom{\frac{1}{2} \bigl(}{}
-|t-s|^{H}\Gamma(1,1)|t-s|^{H^*} \bigr),\nonumber
\end{eqnarray}
where $\Gamma(1,1) = EB_H(1)B_H(1)^*$. Conversely, an OFBM with
covariance function (\ref{e:cov_func_timerevers_OFBM}) is time-reversible.
\end{proposition}
\begin{pf}
This follows from Proposition \ref
{p:time-revers_s.i._processes}(iii).
\end{pf}

\begin{remark}
One should note that for a fixed exponent $H$, not every positive
definite matrix $\Gamma(1,1)$ leads to a valid covariance function
(\ref{e:cov_func_timerevers_OFBM}) for time-reversible OFBMs.

In fact, fix $\Gamma(1,1) = I$, $n=2$. We will show that, for an
exponent of the form
\[
H = \pmatrix{
h & 0\cr
1 & h
},\qquad  h \in(0,1),
\]
there does not exist a time-reversible OFBM $B_{H}$ such that
$EB_{H}(1)B_{H}(1)^* = \Gamma(1,1)$.

From Theorem \ref{t:time-revers_OFBM},
%
\begin{equation}\label{e:Gamma(1,1)}
EB_{H}(1)B_{H}(1)^* = \int_{\bbR} \biggl|\frac{\mathrm{e}^{\mathrm{i}x}-1}{\mathrm{i}x} \biggr|^2 |x|^{-D}
AA^* |x|^{-D^*}\,\mathrm{d}x,
\end{equation}
where $(s_{ij})_{i,j=1,2}:=AA^* \in M(n,\bbR)$.
We have
\begin{eqnarray*}
|x|^{-D} AA^* |x|^{-D^*} &=& |x|^{-2d} \pmatrix{
1 & \cr
\log(x) & 1
}
\pmatrix{
s_{11} & s_{12} \cr
s_{12} & s_{22}
}\pmatrix{
1 & \log(x)\cr
& 1
}
\\
&= &\pmatrix{
s_{11} & s_{11} \log(x) + s_{12} \cr
s_{11} \log(x) + s_{12} & s_{11} (\log(x))^2 + 2 s_{12} \log(x) + s_{22}
}.
\end{eqnarray*}
For notational simplicity, let
\[
r_{k}(d) = \int_{\bbR} \biggl|\frac{\mathrm{e}^{\mathrm{i}x}-1}{\mathrm{i}x} \biggr|^2
(\log(x))^k |x|^{-2d}\,\mathrm{d}x,\qquad  k=0,1,2.
\]
We obtain
\[
EB_{H}(1)B_{H}(1)^* = \pmatrix{
s_{11} r_{0}(d) & s_{11} r_{1}(d) + s_{12} r_{0}(d) \cr
s_{11} r_{1}(d) + s_{12} r_{0}(d) & s_{11} r_{2}(d) + 2 s_{12} r_{1}(d) +
s_{22} r_{0}(d)
}.
\]
On the other hand, for any real symmetric matrix,
the condition for it to have equal eigenvalues is that the
discriminant of the characteristic polynomial is zero. In terms of
$EB_{H}(1)B_{H}(1)^*$,
this means that
\[
s_{11} r_{0}(d) = s_{11} r_{2}(d) + 2 s_{12} r_{1}(d) + s_{22} r_{0}(d),\qquad
s_{11} r_{1}(d) + s_{12} r_{0}(d) = 0.
\]
Therefore, $s_{11} = s_{12} = s_{22} = 0$, which contradicts the
assumption that $\Gamma(1,1) = I$.

This issue is a problem, for instance, in the context of simulation
methods that require
knowledge of the covariance function. For time-reversible OFBMs with
diagonalizable $H$,
one natural way to parameterize $\Gamma(1,1)$ is by means of the
formula (\ref{e:Gamma(1,1)})
since, in this case, the former can be explicitly computed
(see \cite{helgasonpipirasabry2010}).
\end{remark}

Finally, we provide a result (Proposition \ref{p:iso<=>G1=O(n)})
characterizing time-reversibility of some OFBMs in terms of their
symmetry group $G_1$ (see Section \ref{ss:o.s.s._processes}). This
result will be used several times in the next section.

\begin{proposition}\label{p:iso<=>G1=O(n)}
Let $\{B_H(t)\}_{t \in\bbR}$ be an OFBM such that $hI \in
{\mathcal E}(B_H)$ for some $h \in(0,1)$. Then $\{B_H(t)\}_{t
\in\bbR}$ is time-reversible if and only if $G_1(B_H)$ is
conjugate to $O(n)$.
\end{proposition}

\begin{pf}
Regarding necessity, note that if such $B_H$ is time-reversible,
then, by Theorem \ref{t:time-revers_OFBM}, its covariance function
can be written as
\[
\Gamma(t,s) = \int_{\bbR} \biggl(\frac{\mathrm{e}^{\mathrm{i}tx}-1}{\mathrm{i}x}
\biggr) \biggl(\frac{\mathrm{e}^{-\mathrm{i}sx}-1}{-\mathrm{i}x} \biggr)|x|^{-2dI} S \,\mathrm{d}x
\]
for some positive definite $S \in M(n)$ (note that if $S$ is only
positive semidefinite, then the process is not proper).

For sufficiency, consider the covariance function of the OFBM with
exponent $H = hI$, $h \in(0,1)$,
\[
\Gamma(t,s) = \int_{\bbR} \biggl(\frac{\mathrm{e}^{\mathrm{i}tx}-1}{\mathrm{i}x}
\biggr) \biggl(\frac{\mathrm{e}^{-\mathrm{i}sx}-1}{-\mathrm{i}x} \biggr)(x^{-2dI}_{+} AA^* +
x^{-2dI}_{-} \overline{AA^*}) \,\mathrm{d}x.
\]
Define $\widetilde{B}_H = W^{-1}B_H$, where $WO(n)W^{-1} =
G_1(B_H)$ for a positive definite $W$. Then, for any $O \in O(n)$,
\[
\{O\widetilde{B}_H(t)\}_{t \in\bbR}\stackrel{{\mathcal L}}=
\{\widetilde{B}_H(t)\}_{t \in\bbR}.
\]
By the uniqueness of the spectral distribution function, this
implies that $O (W^{-1}AA^*W^{-1} )O^* = W^{-1}AA^*W^{-1}$,
that is, $O(W^{-1}AA^*W^{-1}) = (W^{-1}AA^*W^{-1})O$. Since $O$ is
any matrix in $O(n)$, it follows that $W^{-1}AA^*W^{-1} = cI, c \in
\bbC\backslash\{ 0 \}$ (for a proof of this technical result, see
\cite{didierpipiras2010}). Thus, $AA^* = cW^2$ and $c
> 0$. Hence, $AA^* = \overline{AA^*}$.
\end{pf}

\section{The dichotomy principle}\label{s:dicho}
We now take a closer look at the increments of an OFBM, which form a
stationary process.
\begin{definition}
Let $\{B_{H}(t)\}_{t \in\bbR}$ be an OFBM. The increment process
\[
\{Y_{H}(t)\}_{t \in T} \stackrel{{\mathcal L}}= \{B_{H}(t+1) -
B_{H}(t)\}_{t \in T}, \qquad \mbox{where $T$ = $\bbZ$ or
$\bbR$},
\]
is called \textit{operator fractional Gaussian noise} (OFGN).
\end{definition}

From Theorem \ref{t:spectral_repres_OFBM}, the spectral
representation of OFGN in continuous time is
%
\begin{equation} \label{e:OFGN_spec_representation}
\{Y_{H}(t)\}_{t \in\bbR} \stackrel{{\mathcal L}}=
\biggl\{\int_{\bbR} \mathrm{e}^{\mathrm{i}tx}\frac{\mathrm{e}^{\mathrm{i}x} - 1}{\mathrm{i}x} (x^{-D}_{+}A +
x^{-D}_{-}\overline{A})\widetilde{B}(\mathrm{d}x) \biggr\}_{t \in\bbR}.
\end{equation}
The spectral density of $\{Y_{H}(t)\}_{t \in\bbR}$ is then
%
\begin{equation}
f_{Y_{H}}(x) = \frac{|\mathrm{e}^{\mathrm{i}x}-1|^2}{|x|^2} (x^{-D}_{+}AA^{*}
x^{-D^*}_{+}+ x^{-D}_{-}\overline{AA^{*}}x^{-D^*}_{-}),\qquad  x
\in\bbR,
\end{equation}
since the cross terms are zero.

In discrete time, analogously to the univariate expression,
%
\begin{eqnarray}
&&EY_H(0) Y_H(n)^*\nonumber\\[-8pt]\\[-8pt]
&&\quad  = \int^{\curpi}_{-\curpi} \mathrm{e}^{\mathrm{i}nx} \sum^{\infty}_{k =
-\infty}f_{Y_{H}}(x+2\curpi k) \,\mathrm{d}x,\qquad  n \in\bbZ.\nonumber
\end{eqnarray}
The spectral density of $\{Y_{H}(n)\}_{n \in\bbZ}$ is
then
%
\begin{eqnarray} \label{e:OFGN_spec_dens_disc-time}
g_{Y_{H}}(x) &=& 2 \bigl(1-\cos(x)\bigr)\nonumber\\
&&{}\times\sum^{\infty}_{k =
-\infty} \frac{1}{|x+2\curpi k|^2} \bigl((x+2\curpi k)^{-D}_{+}AA^{*}
(x+2\curpi k)^{-D^*}_{+}\\
&&\hphantom{{}\times\sum^{\infty}_{k =
-\infty} \frac{1}{|x+2\curpi k|^2} \bigl(}
{}+ (x+2\curpi k)^{-D}_{-}\overline{AA^{*}}(x+2\curpi
k)^{-D^*}_{-} \bigr),\qquad  x \in[-\curpi,\curpi].\qquad \nonumber
\end{eqnarray}

The form (\ref{e:OFGN_spec_dens_disc-time}) of the spectral
density leads to the following result.

\begin{theorem} \label{t:dichotomy}
Let $H$ be an exponent with (possibly repeated) characteristic
roots $h_l$, $l=1,\ldots,n$, such that
%
\begin{equation} \label{e:dicho_eigen_between_1/2_and_1}
1/2 < \operatorname{Re}(h_l)<1,\qquad  l=1,\ldots,n.
\end{equation}
Let $g_{Y_{H}}(x) = \{g_{Y_{H}}(x)_{ij}\}$ be the spectral density
(\ref{e:OFGN_spec_dens_disc-time}) of OFGN in discrete time. Then,
for fixed $i,j$, either:
\begin{longlist}[(ii)]
\item[(i)] $|g_{Y_{H}}(x)_{ij}| \rightarrow\infty$, as $x
\rightarrow0$; or
\item[(ii)] $g_{Y_{H}}(x)_{ij} \equiv0$ on
$[-\curpi,\curpi]$.
\end{longlist}
\end{theorem}

\begin{pf}
Let $d_l$ and $N$ be as in (\ref{e:hk_dk}) and (\ref{e:k=1...nN}),
and take $x > 0$. By assumption
(\ref{e:dicho_eigen_between_1/2_and_1}), $0 < \operatorname{Re}(d_l)
< 1/2$. For a given $z > 0$, if we take $-D$ in Jordan canonical
form $PJP^{-1}$, we obtain that
\[
z^{-D} = P \operatorname{diag}
(z^{J_{-d_1}},\ldots,z^{J_{-d_N}}) P^{-1},
\]
where $J_{-d_l}$ is a Jordan block in $J$, $l=1,\ldots, N \leq n$.
Without loss of generality, for $k \geq0$, each term of the
summation (\ref{e:OFGN_spec_dens_disc-time}) involves the matrix
expression
%
\begin{eqnarray}\label{e:Dicho_summ_term}
&&P \operatorname{diag} \bigl((x+2\curpi
k)^{J_{-d_1}}_{+},\ldots,(x+2\curpi k)^{J_{-d_N}}_{+}\bigr)
P^{-1} A
\nonumber\\[-8pt]\\[-8pt]
&&\quad {}\times A^{*}(P^{*})^{-1} \operatorname{diag}\bigl((x+2\curpi
k)^{J^{*}_{-d_1}}_{+},\ldots,(x+2\curpi k)^{J^{*}_{-d_N}}_{+}\bigr)\nonumber
P^{*}.
\end{eqnarray}
Denote the entries of the matrix-valued function
(\ref{e:Dicho_summ_term}) by $h(x+2\curpi k)_{ij}$, $i,j=1,\ldots,n$.
As shown in Appendix \ref{a:Jordan_form}, $h(x+2\curpi k)_{ij}$ is a linear
combination (with complex coefficients) of terms of the form
\[
p(x+2\curpi k)(x+2\curpi k)^{-d_l}q(x+2\curpi k)(x+2\curpi
k)^{-\overline{d}_m},\qquad  l,m = 1,\ldots,N,
\]
where $p(x), q(x)$ are polynomials (with complex coefficients) in
$\log(x)$. Thus,
%
\begin{equation}\label{e:Dicho_summ_finite_k>0}
\sup_{x \in[-\curpi,\curpi]}  \Biggl|\sum^{\infty}_{k=1} \frac{1}{|x+2 \curpi
k|^2} h(x+2\curpi k)_{ij} \Biggr| <
\infty,\qquad  i,j=1,\ldots,n,
\end{equation}
and, therefore,
\begin{eqnarray*}
&&\lim_{x \rightarrow0} 2\bigl(1-\cos(x)\bigr) \Biggl(\frac{1}{|x|^2} h(x)_{ij}
+ \sum^{\infty}_{k=1} \frac{1}{|x+2 \curpi k|^2} h(x+2\curpi
k)_{ij} \Biggr)
\\
&&\quad =\lim_{x \rightarrow0} 2\bigl(1-\cos(x)\bigr) \biggl(\frac{1}{|x|^2}
h(x)_{ij} \biggr).
\end{eqnarray*}
On the other hand, since $\operatorname{Re}(d_l) > 0$ for $l=
1,\ldots,n$, $h(x)_{ij}$ diverges as the power function (times
some $p(x)q(x)$) as $x \rightarrow0$ unless it is identically
zero for all $x$ (in particular, for $x+2 \curpi k$). Thus, by
(\ref{e:Dicho_summ_finite_k>0}) and the fact that
$\frac{2(1-\cos(x))}{x^2} \rightarrow1$ as $x \rightarrow0$, the
claim
follows.
\end{pf}

In the univariate context, the range $(1/2,1)$ for $H$ is commonly
known as that of long-range dependence (LRD). In the multivariate
context, characteristic roots of $H$ with real parts between $1/2$
and 1 have the potential to generate divergence of the spectrum at
zero. Theorem \ref{t:dichotomy} thus states that if OFGN is long-range
dependent in the sense of
(\ref{e:dicho_eigen_between_1/2_and_1}), then the cross
correlation between any two components is characterized by the
following \textit{dichotomy}:
\begin{itemize}
\item it either has a divergent spectrum at zero, a characteristic
usually associated with LRD; or
\item it is identically equal to
zero.
\end{itemize}

Obtaining a similar dichotomy principle for a larger range of
characteristic roots than that in
(\ref{e:dicho_eigen_between_1/2_and_1}) is much more delicate. The
following two examples illustrate some of the potential
difficulties. Example \ref{ex:dicho_BM} shows that if one of the
characteristic roots of $D=H-(1/2)I$ is 0, then the dichotomy may not
hold. Example \ref{ex:dicho_d_-d} shows that certain cancellations
may occur in the cross spectrum if the characteristic roots of $D$
have opposite signs.

\begin{example}\label{ex:dicho_BM}
If, for instance, $D = P\operatorname{diag}(d,0)P^{-1}$, where $0 <
d <1/2$,
\[
P= \pmatrix{
1 & \sqrt{2}/2 \cr
0 & \sqrt{2}/2 }
\]
and $A := P$, then, as $x \rightarrow0$,
\[
g_{Y_H}(x) \sim \pmatrix{
x^{-2d}+1/2 & 1/2 \cr
1/2 & 1/2 },
\]
where $\sim$ indicates entrywise asymptotic equivalence. As a
consequence, if one of the components of OFBM behaves like
Brownian motion, then this may create cross short-range dependence
among the components. This example is a direct consequence of a
more general operator parameter $D$ whose eigenspaces are not the
canonical axes. If we take, instead, $D = \operatorname{diag}(d,0)$,
whose eigenspaces are the canonical axes, each term of the
summation (\ref{e:OFGN_spec_dens_disc-time}) has the form
\[
\operatorname{diag} \bigl((x+2\curpi k)^{-d}_{\pm},0\bigr) AA^*
\operatorname{diag}\bigl((x+2\curpi k)^{-d}_{\pm},0\bigr)
= \pmatrix{ s_{11}(x+2\curpi k)^{-2d}_{\pm} & 0\cr
0 & 0},
\]
where
%
\begin{equation}\label{e:notation_sij=AA*}
(s_{ij})_{i,j=1,2} := AA^*
\end{equation}
and thus the dichotomy holds.
\end{example}

\begin{example}\label{ex:dicho_d_-d}
Consider $A \in \mathit{GL}(n,\bbR)$ and $D = \operatorname{diag}(d,-d)$,
where $d \in(0,1/2)$. Using the notation
(\ref{e:notation_sij=AA*}), we have
\[
g_{Y_H}(x) \sim P \operatorname{diag}(x^{-d},x^{d})P^* AA^* P
\operatorname{diag}(x^{-d},x^{d})P^* = \pmatrix{
s_{11} x^{-2d} & s_{12} \cr
s_{12} & s_{22} x^{2d}
}
\]
as $x \rightarrow0$. Here, the multivariate differencing effects
of the operator $D$ cancel out in the cross-entries.
\end{example}

From a practical perspective, Theorem \ref{t:dichotomy} raises the
question of whether
the class of OFGNs is flexible enough to capture multivariate LRD structures.
This, and related issues regarding multivariate discrete time series,
will be explored in future work.


\section{Operator Brownian motions}\label{s:OBM}

In this section, we shall adopt the following definition of
multivariate Brownian motion and establish some of its
properties.
\begin{definition}\label{def:OBM}
The proper process $\{B_H(t)\}_{t \in
\bbR}$ is an operator Brownian motion (OBM) if it is a Gaussian
o.s.s.\ process which has stationary and independent increments and
satisfies $B_H(0)=0$ a.s.
\end{definition}

In place of the condition $B_H(0)=0$ a.s., we can assume that the
characteristic roots of the o.s.s.\ exponent
$H$ have positive real parts, which implies the former condition.
Another important way to motivate Definition \ref{def:OBM} is as
follows: since $\{B_H(t)\}_{t \in
\bbR}$ is $L^2$-continuous (see the beginning of the proof of Theorem
\ref{t:spectral_repres_OFBM}), by Hudson and Mason \cite
{hudsonmason1981JrMultAnalysis}, Theorem 4, and Hudson and Mason
\cite{hudsonmason1982}, Theorem 7, our Definition \ref{def:OBM}
implies that $(1/2)I$ can always be taken
as an  exponent of OBM.

The next proposition and example show that an OFBM $B_H$ with $(1/2)I
\in{\mathcal E}(B_H)$ is not necessarily
an OBM. This stands in contrast with the univariate case.

\begin{proposition}\label{p:OBM<=>MN*=NM*}
Let $\{B_H(t)\}_{t \in\bbR}$ be an OFBM with exponent $H=(1/2)I$.
Consider its time domain
representation (\ref{e:time domain_OFBM_D=0}) with parameters $M$ and
$N$, or its spectral domain representation~(\ref{e:spectral_repres_OFBM}) with $A=A_1 + \mathrm{i} A_2$. Then $\{B_H(t)\}_{t \in
\bbR}$ is an OBM if and only if the following two equivalent conditions hold:
\begin{longlist}[(ii)]
\item[(i)] $MN^{*}=NM^{*}$;
\item[(ii)] $A_2 A^*_1 = A_1 A^*_2$.
\end{longlist}
\end{proposition}
\begin{pf}
Since $(1/2)I \in{\mathcal E}(B_H)$, it follows that $B_H(0)=0$ a.s.
Therefore, we only have to establish that the increments are independent
if and only if (i) holds. Demonstrating the equivalence between (i)
and (ii) is straightforward by using the relation (\ref{e:M_N}).

Write the time domain representation (\ref{e:time domain_OFBM_D=0}) of
$B_H$ as
%
\begin{equation}\label{e:OBM<=>MN*=NM*}
\int_{\bbR} \bigl(2\bigl(1_{\{t-u>0\}}-1_{\{-u>0\}}\bigr)M +
(\log|t-u|-\log|-u|)N \bigr)B(\mathrm{d}u).
\end{equation}
Take $s_1 < t_1 < s_2 < t_2$. For the increments of the process $B_H$,
we have
%
\begin{eqnarray}\label{e:integ_OBM}
&&E\bigl(B_H(t_1)-B_H(s_1)\bigr)\bigl(B_H(t_2)-B_H(s_2)\bigr)^*\nonumber
\\
&&\quad = \int_{\bbR} \biggl(4 \cdot1_{\{s_1 < u < t_1\}}1_{\{s_2 < u < t_2\}}MM^*
+ \log\frac{|t_1-u|}{|s_1-u|}\log\frac{|t_2-u|}{|s_2-u|}NN^*
\\
&&\hphantom{\int_{\bbR} \biggl(}\qquad {}+ 2\cdot 1_{\{s_1 < u < t_1\}}\log\frac{|t_2-u|}{|s_2-u|}MN^* + 2\cdot 1_{\{s_2
< u < t_2\}}\log\frac{|t_1-u|}{|s_1-u|}NM^* \biggr) \,\mathrm{d}u.\nonumber
\end{eqnarray}
From the univariate time domain representation of Brownian motion,
the first two of the four terms in (\ref{e:integ_OBM}) have zero
integral. Define $\varphi(u) = u (\log(u)-1)$. The right-hand side
of the expression (\ref{e:integ_OBM}) then becomes
\begin{eqnarray*}
&&\bigl(\varphi(t_2-s_1)
-\varphi(t_2-t_1)-\varphi(s_2-s_1)+\varphi(s_2-t_1)\bigr)MN^*
\\
&&\quad {}+ \bigl(\varphi(t_2-t_1)
-\varphi(s_2-t_1)-\varphi(t_2-s_1)+\varphi(s_2-s_1)\bigr)NM^*,
\end{eqnarray*}
which is identically zero if and only if $MN^{*}=NM^{*}$.
\end{pf}

\begin{example}\label{ex:notOBM,not_timerevers}
Consider a bivariate process $X$ defined by the expression
(\ref{e:OBM<=>MN*=NM*}). Set $M = I$ and let $N=L \in \mathit{so}(2)
\backslash\{0\}$. Then $MN^* = - NM^* \neq0$ (from which the
cross terms in expression (\ref{e:integ_OBM}) cancel out when $s_1
= s_2 = 0$ and $t_1 = t_2 = t$) and
\[
EX(t)X(t)^* = 4 |t|I + \curpi^2 |t| L(-L) = |t| (4I-\curpi^2 L^2),
\]
which is a full rank matrix for $t \neq0$ (to obtain the constant
$\curpi^2$, one can use, e.g., Proposition~9.2 in
\cite{taqqu2003} and Proposition 5.1 in
\cite{stoevtaqqu2007}). Hence, $X$ is proper. This gives an
example of an OFBM for which $(1/2)I \in{\mathcal E}(X)$ but
which is \textit{not} an OBM. Moreover, it is an example of an
OFBM with an exponent of the form $hI$, $h \in(0,1)$, but which
is \textit{not} time-reversible and for which $G_1 \cong O(2)$
does \textit{not} hold by Proposition \ref{p:iso<=>G1=O(n)} (in
contrast, by Hudson and Mason
\cite{hudsonmason1982}, Theorem~6, $G_1(X) \cong O(n)$ implies that
$hI \in
{\mathcal E}(X)$ for some $h$).
%
\end{example}

The following result is a direct consequence of Theorem
\ref{t:time-revers_OFBM} and Proposition \ref{p:OBM<=>MN*=NM*}. It
shows that in the class of OFBMs with exponent $H=(1/2)I$,
time-reversibility is equivalent to independence of increments.

\begin{corollary}\label{c:H=1/2_timerevers<=>OBM}
Let $\{B_H(t)\}_{t \in\bbR}$ be a time-reversible OFBM with exponent
$H=(1/2)I$.
Then $\{B_H(t)\}_{t \in\bbR}$ is an OBM. Conversely, let $\{B_H(t)\}
_{t \in\bbR}$ be an OBM. It is then time-reversible (and
has exponent $H=(1/2)I$).
\end{corollary}

\begin{remark}\label{r:equivalence_timerevers_G1=O(n)}
Note that, as a consequence of Proposition \ref{p:iso<=>G1=O(n)},
time-reversibility may be replaced in Corollary
\ref{c:H=1/2_timerevers<=>OBM} by the condition that $G_1$ is
conjugate to $O(n)$. In other words, an OBM is elliptically
symmetric.
\end{remark}

We conclude by providing a spectral representation for OBM.

\begin{proposition}\label{p:spectral_repres_OBM}
Let $\{B_H(t)\}_{t \in\bbR}$ be an OBM. Then
%
\begin{equation}\label{e:spectral_repres_OBM}
\{B_H(t)\}_{t \in\bbR} \stackrel{{\mathcal L}}=
\biggl\{\int_{\bbR}\frac{\mathrm{e}^{\mathrm{i}tx}-1}{\mathrm{i}x} W\widetilde{B}(\mathrm{d}x) \biggr\}_{t
\in\bbR} \stackrel{{\mathcal L}}= \{W B(t)\}_{t \in\bbR}
\end{equation}
for some positive definite operator $W$, where $\{B(t)\}_{t \in
\bbR}$ is a vector of independent standard BMs.
\end{proposition}
\begin{pf}
Consider the spectral domain representation of $B_H$ with parameter $A
= A_1 + \mathrm{i} A_2$. Set $W:= (A_1 A^*_1 + A_2 A^*_2)^{1/2}$. The result follows
from Proposition \ref{p:OBM<=>MN*=NM*}(ii), and relations (\ref
{e:cov_struct_BH(t)}) and (\ref{e:AA*}).\vspace*{-1pt}
\end{pf}

\begin{appendix}

\section{Fourier transforms of OFBM kernels}\label{ss:FT}

In this appendix, the goal is to prove Proposition
\ref{p:FT_timedom_filter}. First, we state a condensed version of
Horn and Johnson
\cite{hornjohnson1991}, Theorems 6.2.9 and 6.2.10, pages 412--416,
which will be useful in
the subsequent derivations. We shall use the notation introduced before
Definition \ref{d:matrix_function}.\vspace*{-1pt}
\begin{theorem}\label{t:Theo6.2.9_and_6.2.10_in_HornJohnson}
Let $f, g\dvtx U \rightarrow\bbC$ be two stem functions and let ${\mathcal
M}_{fg}= {\mathcal M}_f \cap{\mathcal M}_g$. Then:
\begin{enumerate}[(iii)]
\item[(i)] the primary matrix function $f\dvtx{\mathcal M}_f \rightarrow
M(n,\bbC)$ is well defined in the sense that the value of $f(\Lambda)$,
$\Lambda\in{\mathcal M}_f $, is independent of the particular Jordan
canonical form (i.e., block permutation) used to represent it;
\item[(ii)] $f(\Lambda) = g(\Lambda)$ if and only if $f^{(j)}(\lambda
_k) = g^{(j)}(\lambda_k)$ for $j = 0,1,\dots,r_k - 1$, $k = 1,\dots,N$
and $\Lambda\in{\mathcal M}_{fg}$;
\item[(iii)] for $q(z):=f(z)g(z)$, we have
$q(\Lambda)=f(\Lambda)g(\Lambda)=g(\Lambda)f(\Lambda)$ for $\Lambda\in
{\mathcal M}_{fg}$;
\item[(iv)] for $s(z):=f(z)+g(z)$, we have $s(\Lambda)=f(\Lambda)+g(\Lambda
)$ for $\Lambda\in{\mathcal M}_{fg}$.\vspace*{-1pt}
\end{enumerate}
\end{theorem}

Throughout this section, we assume (\ref{e:reg_cond_char_roots_H})
and (\ref{e:ch_roots_neq_1/2}). Denote by ${\mathcal F}$ the
Fourier transform operator. For $d \in\bbC$ such that
%
\begin{equation}\label{e:range_d}
\operatorname{Re}(d) \in(-1/2,1/2) \backslash\{0\},\vspace*{-1pt}
\end{equation}
define
\[
f_{\pm}(t,u,d) = (t-u)^{d}_{\pm}-(-u)^{d}_{\pm}\vspace*{-1pt}
\]
and
\[
h_{\pm}(t,x,d)= \frac{\mathrm{e}^{\mathrm{i}tx}-1}{\mathrm{i}x}|x|^{-d}\Gamma(d+1) \mathrm{e}^{\mp
\operatorname{sign}(x){\mathrm{i} \curpi d}/{2}}.\vspace*{-1pt}
\]
It is well known that
%
\begin{equation} \label{e:Tinv_timedom_filter}
{\mathcal F}(f_{\pm}(t,\cdot,d))(x) =
h_{\pm}(t,x,d)\vadjust{\goodbreak}
\end{equation}
when $d \in(-1/2, 1/2) \backslash\{0\}$ (see, e.g.,
\cite{pipirastaqqu2003}, page\ 175). One can
show that (\ref{e:Tinv_timedom_filter}) also holds under
(\ref{e:range_d}) (see Remark \ref{r:on_analytic_continuation}).

For the purpose of calculating Fourier transforms of primary
matrix functions associated with the stem functions $f_{\pm}$ and
$h_{\pm}$, we will need to consider derivatives of the latter with
respect to~$d$. Note that, for fixed $x$, the functions $\Gamma(d
+1)$, $\mathrm{e}^{\mp\operatorname{sign}(x){\mathrm{i} \curpi d}/{2}}$ and $|x|^{- d}$
are holomorphic on the domain $-\frac{1}{2} < \operatorname{Re}(d) <
\frac{1}{2}$. Thus, so are the functions $h_{\pm}(t,x,d)$. Note
that, for fixed $t$ and~$u$, since $(t-u)^{d}_{\pm}$ and
$(-u)^{d}_{\pm}$ are holomorphic on the domain $-1/2 <
\operatorname{Re}(d)<1/2$, then so are $f_{\pm}(t,u,d)$.

As a consequence, by Theorem \ref{t:Theo6.2.9_and_6.2.10_in_HornJohnson}(i)--(iv),
\[
h_{\pm}(t,x,D) = \frac{\mathrm{e}^{\mathrm{i}tx}-1}{\mathrm{i}x}|x|^{-D}\Gamma(D+I)\mathrm{e}^{\mp
\operatorname{sign}(x) \mathrm{i} \curpi D/{2}}
\]
and
\[
f_{\pm}(t,u,D) = (t-u)^{D}_{\pm}-(-u)^{D}_{\pm}.
\]
We now need to show that
%
\begin{equation}\label{e:g1or2=IFT(timedom_filters)}
{\mathcal F}(f_{\pm}(t,u,D))(x)= h_{\pm}(t,x,D),
\end{equation}
where ${\mathcal F}$ is the \textit{entrywise} Fourier transform
operator.

\begin{pf*}{Proof of Proposition \ref{p:FT_timedom_filter}}
We will break up the proof into three cases:

\textit{Case} 1: $-1/2<\operatorname{Re}(d_k)< 0$, $k = 1,\ldots,N$. We will
develop the calculations for $h_+$, which can be easily adapted to
$h_{-}$. By Theorem \ref{t:Theo6.2.9_and_6.2.10_in_HornJohnson}(ii), in
the case of $h_{+}$,
(\ref{e:g1or2=IFT(timedom_filters)}) is equivalent to
%
\begin{equation}\label{e:diff_under_integ}
\frac{\partial^{j}}{\partial d^j }h_{+}(t,x,d) = \frac{\partial
^{j}}{\partial d^j }\int_{\bbR}
\mathrm{e}^{\mathrm{i}ux} f_{+}(t,u,d)\,\mathrm{d}u = \int_{\bbR}
\mathrm{e}^{\mathrm{i}ux} \frac{\partial^{j}}{\partial d^j }f_{+}(t,u,d)\,\mathrm{d}u
\end{equation}
at $d = d_k$, for $k=1,\ldots,N$, $j = 0,1,\ldots,r_k-1$. Consider the
domain $\Delta(\underline{d},\overline{d}) := \{d \in\bbC\dvtx
\underline{d} < \operatorname{Re}(d) < \overline{d})\}$, where $-1/2
< \underline{d} < \overline{d} < 0$, which is open and convex.
Consider $j =1$, that is, the first derivative. Fix $d^* \in
\Delta(\underline{d},\overline{d})$ and take a sequence
$\{d_{m}\}_{m \in\bbN} \subseteq
\Delta(\underline{d},\overline{d})$ such that $d_{m} \rightarrow
d^*$. For each $m$, by the mean value theorem for holomorphic
functions (\cite{evardjafari1992}, Theorem 2.2), there exist constants $\delta_i(m)$, where
\[
\delta_i(m) = \alpha_{m,i} d_m + (1 - \alpha_{m,i}) d^*,\qquad
\alpha_{m,i} \in(0,1), i=1,2,
\]
such that
\begin{eqnarray*}
\frac{f_{+}(t,u,d_m)-f_{+}(t,u,d^*)}{d_m - d^*} &=&
\operatorname{Re} \biggl(\frac{f_{+}(t,u,d_m)-f_{+}(t,u,d^*)}{d_m -
d^*} \biggr)\\
&&{} + \mathrm{i}
\operatorname{Im} \biggl(\frac{f_{+}(t,u,d_m)-f_{+}(t,u,d^*)}{d_m -
d^*} \biggr)
\\
&=& \operatorname{Re} \biggl(\frac{\partial}{\partial
d}f_{+}(t,u,\delta_1(m)) \biggr) + \mathrm{i}
\operatorname{Im} \biggl(\frac{\partial}{\partial
d}f_{+}(t,u,\delta_2(m)) \biggr).
\end{eqnarray*}

We will now obtain an integrable function that majorizes
$|\frac{\partial}{\partial d}f_{+}(t,\cdot,d)|$ for all $d \in
\Delta(\underline{d},\overline{d})$. Assume, without loss of
generality, that $t
> 0$, and take $d \in\Delta(\underline{d},\overline{d})$ and $\delta
> 0$ such that $-1/2 < \operatorname{Re}(\underline{d}) - \delta$
and $\operatorname{Re}(\overline{d}) + \delta< 0$. From the continuity of
$\frac{\partial}{\partial d}f_{+}(t,u,d)$
for $0 \leq u < t$, there exist constants $K_1$ and $\eta_{1}$
such that
\begin{eqnarray} \label{e:deriv_f+_bound1}
\biggl| \frac{\partial}{\partial d}f_{+}(t,u,d) \biggr| &\leq&
|\log(t-u)_+||(t-u)^{d}_{+}|
\nonumber\\[-8pt]\\[-8pt]
&\leq& K_1 1_{[0,t-\eta_{1}]}(u) +
\bigl|(t-u)^{\operatorname{Re}(\underline{d})-\delta}_{+} \bigr|1_{(t-\eta_{1},t)}(u),\qquad
u \geq0.\nonumber
\end{eqnarray}
Also, there exists a constant $K_2$ such that
%
\begin{equation}\label{e:deriv_f+_bound2}
|\log(t-u)_+ (t-u)^{d}_{+} - \log(-u)_+ (-u)^{d}_{+}| \leq K_2
+\bigl|(-u)^{\operatorname{Re}(\underline{d}) - \delta}_{+}\bigr|,\qquad  -1
\leq u < 0.
\end{equation}
One can show that there exist constants $K_3$ and $\eta_2 < -1$
such that
%
\begin{equation}\label{e:deriv_f+_bound3}
|\log(t-u)_+ (t-u)^{d}_{+} - \log(-u)_+ (-u)^{d}_{+}| \leq K_3
(-u)^{\operatorname{Re}(\overline{d}) + \delta-1}_{+},\qquad  u <
\eta_2.
\end{equation}
From (\ref{e:deriv_f+_bound1}), (\ref{e:deriv_f+_bound2}) and
(\ref{e:deriv_f+_bound3}), and from the fact that
$\frac{\partial}{\partial d}f_{+}(t,u,d)$ is bounded on $\eta_2
\leq u \leq-1$ uniformly in $d \in
\Delta(\underline{d},\overline{d})$, we conclude that the ratio
$\frac{f_{+}(t,\cdot,d_m)-f_{+}(t,\cdot,d^*)}{d_m - d^*}$ is
bounded by a function in $L^1(\bbR)$. Thus, by the dominated
convergence theorem (for $\bbC$-valued functions),
we have
\[
\int_{\bbR}\mathrm{e}^{\mathrm{i}ux}\frac{f_{+}(t,u,d_m)-f_{+}(t,u,d^*)}{d_m -
d^*} \,\mathrm{d}u \rightarrow\int_{\bbR}\mathrm{e}^{\mathrm{i}ux}\frac{\partial}{\partial
d}f_{+}(t,u,d^*) \,\mathrm{d}u,\qquad  m \rightarrow\infty.
\]
We can always assume that $t \neq0$ and the case of $t < 0$ can
be dealt with in a similar fashion. The extension of the above
argument for derivatives of higher order $j$ poses no additional
technical difficulties. This establishes
(\ref{e:diff_under_integ}).

\textit{Case} 2: $0<\operatorname{Re}(d_k)< 1/2$, $k = 1,\ldots,N$. In this
range, the upper bound in (\ref{e:deriv_f+_bound3}) is not in
$L^1(\bbR)$ so we need a slightly different procedure. Since
$h_{\pm}(t,\cdot,d) \in L^2(\bbR)$, we can apply ${\mathcal
F}^{-1}$ on both sides of equation (\ref{e:Tinv_timedom_filter})
and obtain
\[
f_{\pm}(t,u,d) = {\mathcal F}^{-1}(h_{\pm}(t,\cdot,d))(u).
\]
Therefore, it suffices to show that
%
\begin{equation}\label{e:q(t,x,D)=IFT(h(t,x,D))}
f_{\pm}(t,u,D) = {\mathcal F}^{-1}(h_{\pm}(t,x,D)),
\end{equation}
where ${\mathcal F}^{-1}$ is the \textit{entrywise} inverse
Fourier transform.

Note that expression (\ref{e:q(t,x,D)=IFT(h(t,x,D))}) is equivalent to
\[
\frac{\partial^j}{\partial d^j}f_{\pm}(t,x,d) = \frac{\partial
^j}{\partial d^j}\int_{\bbR}\mathrm{e}^{-\mathrm{i}ux}
h_{\pm}(t,x,d) \,\mathrm{d}x = \int_{\bbR}\mathrm{e}^{-\mathrm{i}ux}
\frac{\partial^j}{\partial d^j}h_{\pm}(t,x,d) \,\mathrm{d}x
\]
at $d = d_k$, for $k=1,\ldots,N$, $j = 0,1,\ldots,r_k-1$. To show
this, one may proceed as in the case of $-1/2 < \operatorname{Re}(d)
< 0$. The existence of an upper bound in $L^1(\bbR)$ is ensured by
the fact that
%
\begin{eqnarray} \label{e:upperbound_derivative}
\biggl| \frac{\partial^j}{\partial
d^j} \biggl( \biggl(\frac{\mathrm{e}^{\mathrm{i}tx}-1}{\mathrm{i}x} \biggr)|x|^{-d} \biggr) \biggr|
&\leq&
\biggl| \biggl(\frac{\mathrm{e}^{\mathrm{i}tx}-1}{\mathrm{i}x} \biggr)|\log|x||^j
|x|^{-\operatorname{Re}(d)} \biggr|\nonumber
\\
&\leq& \biggl| \biggl(\frac{\mathrm{e}^{\mathrm{i}tx}-1}{\mathrm{i}x} \biggr)|\log|x||^{j}
|x|^{-\operatorname{Re}(\overline{d})} 1_{\{0 < |x| \leq1\}} \biggr|\\
&&{} +
\biggl| \biggl(\frac{\mathrm{e}^{\mathrm{i}tx}-1}{\mathrm{i}x} \biggr)|\log|x||^j
|x|^{-\operatorname{Re}(\underline{d})} 1_{\{1 < |x| <
\infty\}} \biggr|,\nonumber
\end{eqnarray}
which is integrable for all $d \in\Delta(\underline{d},
\overline{d})$, $0 < \underline{d} < \overline{d} < 1/2$.\\

\textit{General case}: As a consequence of
(\ref{e:g1or2=IFT(timedom_filters)}),
%
\begin{equation}\label{e:IFT_f_pm}
{\mathcal F}(f_{\pm}(t,\cdot,J))(x) = h_{\pm}(t,x,J)
\end{equation}
holds, where $J$ is a matrix in Jordan canonical form with
characteristic roots satisfying (\ref{e:reg_cond_char_roots_H})
and (\ref{e:ch_roots_neq_1/2}). Now pre- and post-multiply
equation (\ref{e:IFT_f_pm}) by $P$ and $P^{-1}$, respectively. Since
\[
P\Gamma(-J)P^{-1} = \Gamma(-D),\qquad
P\mathrm{e}^{(\mathrm{i}\curpi/2)(J+I)}P^{-1} = \mathrm{e}^{(\mathrm{i}\curpi/{2})(D+I)},
\]
it follows that
\[
Pf_{\pm}(t,u,J)P^{-1} = f_{\pm}(t,u,D),\qquad
Ph_{\pm}(t,u,J)P^{-1} = h_{\pm}(t,u,D),
\]
from which we obtain equation
(\ref{e:g1or2=IFT(timedom_filters)}).
\end{pf*}

\begin{remark}\label{r:on_analytic_continuation}
A common way to prove that (\ref{e:Tinv_timedom_filter}) also
holds for $d \in\bbC$ satisfying (\ref{e:range_d}) is by analytic
continuation. In particular, this requires the ability to
differentiate under the integral sign in the Fourier transform.
The latter could be achieved by following the argument in the
proof of Proposition \ref{p:FT_timedom_filter}.
\end{remark}

\section{Some useful integrals}\label{a:useful_integ}

In this appendix, we calculate the inverse Fourier transforms used
in the proof of Theorem \ref{t:time domain_OFBM}(ii) and
Example \ref{ex:time domain_D=Jordan_and_eig=0}. We shall use
several formulas from \cite{gradshteynryzhik2007}:
\begin{eqnarray}\label{e:int_indicator_sin}
&\displaystyle\int_{\bbR} 1_{\{x<0\}}\frac{\sin(ax)}{x} \,\mathrm{d}x = \int_{\bbR}
1_{\{x>0\}}\frac{\sin(ax)}{x} \,\mathrm{d}x =
\frac{\curpi}{2}\operatorname{sign}(a)\qquad  (\mbox{page\ 423}),
\\
\label{e:int_indicator>0_cos}
&\displaystyle\int_{\bbR} 1_{\{x>0\}}\frac{\cos(ax)-\cos(bx)}{x} \,\mathrm{d}x =
\log\frac{|b|}{|a|}\qquad  (\mbox{page\ 447})
\\
\label{e:int_indicator<0_cos}
&\displaystyle\biggl(\mbox{therefore, }\int_{\bbR} 1_{\{x<0\}}\frac{\cos(ax)-\cos(bx)}{x} \,\mathrm{d}x =
-\log\frac{|b|}{|a|} \biggr),
\\
\label{e:int_log_sin}
&\displaystyle\int^\infty_0 \log(x) \sin(ax) \frac{\mathrm{d}x}{x} = - \frac{\curpi}{2}\bigl(C +
\log(a)\bigr),\qquad  a > 0\qquad  (\mbox{page\ 594})
\\
\label{e:int_log_sin_any_a}
&\displaystyle\biggl(\mbox{therefore, with $a \in\bbR$, }\int^\infty_0 \log(x) \sin(ax) \frac{\mathrm{d}x}{x} = \int^{0}_{-\infty}
\log(x_{-}) \sin(ax) \frac{\mathrm{d}x}{x}\nonumber\\[-8pt]\\[-8pt]
&\displaystyle\hspace*{200pt} = - \frac{\curpi}{2}(C +
\log|a|)\operatorname{sign}(a) \biggr),\nonumber
\\\label{e:int_log_cos}
&\displaystyle\int^\infty_0 \log(x) \bigl(\cos(ax)-\cos(bx)\bigr) \frac{\mathrm{d}x}{x} =
\log\biggl(\frac{a}{b} \biggr) \biggl(C + \frac{1}{2}\log(ab) \biggr),\qquad
a, b > 0\nonumber\\[-8pt]\\[-8pt]
&\hspace*{-135pt}\qquad  (\mbox{page\ 594}),\hspace*{135pt}\nonumber
\end{eqnarray}
where $C$ is Euler's constant
%
\begin{equation}\label{e:int_log_cos_<0}
\biggl(\mbox{therefore, }\int^{0}_{-\infty} \log(x_{-}) \bigl(\cos(ax)-\cos(bx)\bigr) \frac{\mathrm{d}x}{x} =
-\log\biggl(\frac{|a|}{|b|} \biggr) \biggl(C + \frac{1}{2}\log(|ab|) \biggr)\biggr).
\end{equation}

Using these formulas, we obtain that, for $x > 0$,
\begin{eqnarray*}
&&\int_{\bbR}\mathrm{e}^{-\mathrm{i}ux}\frac{\mathrm{e}^{\mathrm{i}tx}-1}{\mathrm{i}x}1_{\{x>0\}}\,\mathrm{d}x
\\
&&\quad = \int_{\bbR}\frac{1}{\mathrm{i}x} \biggl(\cos\bigl((t-u)x\bigr)-\cos(ux)+\mathrm{i}\bigl(\sin\bigl((t-u)x\bigr)+\sin
(ux)\bigr) \biggr)1_{\{x>0\}}\,\mathrm{d}x
\\
&&\quad = \frac{1}{\mathrm{i}}\log\biggl(\frac{|u|}{|t-u|} \biggr) +
\frac{\curpi}{2} \bigl(\operatorname{sign}(t-u)-
\operatorname{sign}(-u)\bigr).
\end{eqnarray*}
Similarly, for $x < 0$,
\begin{eqnarray*}
&&\int_{\bbR}\mathrm{e}^{-\mathrm{i}ux}\frac{\mathrm{e}^{\mathrm{i}tx}-1}{\mathrm{i}x}1_{\{x<0\}}\,\mathrm{d}x
\\
&&\quad = \int_{\bbR}\frac{1}{\mathrm{i}x} \biggl(\cos\bigl((t-u)x\bigr)-\cos(ux)+\mathrm{i}\bigl(\sin\bigl((t-u)x\bigr)+\sin
(ux)\bigr) \biggr) 1_{\{x<0\}}\,\mathrm{d}x
\\
&&\quad =-\frac{1}{\mathrm{i}}\log\biggl(\frac{|u|}{|t-u|} \biggr) +
\frac{\curpi}{2} \bigl(\operatorname{sign}(t-u)-\operatorname{sign}(-u)\bigr).
\end{eqnarray*}
Therefore,
\begin{eqnarray} \label{e:IFT_spec_filter_d=0}
&&\int_{\bbR}\mathrm{e}^{-\mathrm{i}ux}\frac{\mathrm{e}^{\mathrm{i}tx}-1}{\mathrm{i}x}\bigl(1_{\{x>0\}}A +
1_{\{x<0\}}\overline{A}\bigr) \,\mathrm{d}x
\nonumber\\[-8pt]\\[-8pt]
&&\quad =\bigl(\operatorname{sign}(t-u)-
\operatorname{sign}(-u)\bigr)\frac{1}{2} \operatorname{Re}(A)+\log\biggl(\frac{|u|}{|t-u|} \biggr)
\frac{1}{\curpi}\operatorname{Im}(A),\nonumber
\end{eqnarray}
which is the formula used in the proof of Theorem \ref{t:time domain_OFBM}(ii).

We now turn to the calculations of the inverse Fourier transform
(\ref{e:IFT_D=Jordan_and_eig=0}) in Example \ref{ex:time
domain_D=Jordan_and_eig=0}. Note that
\[
|x|^{-D}= \pmatrix{ 1 & 0\cr
-\log|x| & 1},\qquad  x > 0.
\]
We only need to calculate the inverse Fourier transform of the log
term on the lower off-diagonal.

For $x > 0$, using the formulas above,
\begin{eqnarray*}
&&\int_{\bbR}\mathrm{e}^{-\mathrm{i}ux}\frac{\mathrm{e}^{\mathrm{i}tx}-1}{\mathrm{i}x}\log(x_+) 1_{\{x>0\}}\,\mathrm{d}x
\\
&&\quad = \int_{\bbR}\frac{1}{\mathrm{i}x} \bigl(\cos\bigl((t-u)x\bigr)-\cos(ux)+\mathrm{i}\bigl(\sin\bigl((t-u)x\bigr)+\sin
(ux)\bigr) \bigr)\log(x_+) 1_{\{x>0\}}\,\mathrm{d}x
\\
&&\quad =\frac{1}{\mathrm{i}}\log\biggl(\frac{|t-u|}{|u|} \biggr) \biggl(C +
\frac{1}{2}\log(|t-u||u|) \biggr)\\
&&\qquad {} -
\frac{\curpi}{2}\bigl((C+\log|t-u|)\operatorname{sign}(t-u) -
(C+\log|u|)\operatorname{sign}(-u)\bigr).
\end{eqnarray*}
Similarly, for $x < 0$,
\begin{eqnarray*}
&&\int_{\bbR}\mathrm{e}^{-\mathrm{i}ux}\frac{\mathrm{e}^{\mathrm{i}tx}-1}{\mathrm{i}x}\log(x_-) 1_{\{x<0\}}\,\mathrm{d}x
\\
&&\quad = \int_{\bbR}\frac{1}{\mathrm{i}x} \bigl(\cos\bigl((t-u)x\bigr)-\cos(ux)+\mathrm{i}\bigl(\sin\bigl((t-u)x\bigr)+\sin
(ux)\bigr) \bigr)\log(x_-) 1_{\{x<0\}}\,\mathrm{d}x
\\
&&\quad = - \frac{1}{\mathrm{i}}\log\biggl(\frac{|t-u|}{|u|} \biggr) \biggl(C +
\frac{1}{2}\log(|t-u||u|) \biggr)\\
&&\qquad {} -
\frac{\curpi}{2}\bigl((C+\log|t-u|)\operatorname{sign}(t-u) -
(C+\log|u|)\operatorname{sign}(-u)\bigr).
\end{eqnarray*}
Therefore, for $a \in\bbC$,
%
\begin{eqnarray}\label{e:IFT_spec_filter_d=0_log_on_the_offdiag}
&&\int_{\bbR}\mathrm{e}^{-\mathrm{i}ux}\frac{\mathrm{e}^{\mathrm{i}tx}-1}{\mathrm{i}x}\bigl(-\log(x_+)
1_{\{x>0\}}a - \log(x_-) 1_{\{x<0\}}\overline{a}\bigr)
\,\mathrm{d}x\nonumber
\\
&&\quad = \bigl((C+{\log}|t-u|)\operatorname{sign}(t-u) -
(C+{\log}|u|)\operatorname{sign}(-u) \bigr)
\frac{1}{2}\operatorname{Re}(a)
\\
&&\qquad {}+ \log\biggl( \frac{|t-u|}{|u|} \biggr) \biggl(C +
\frac{1}{2}\log(|t-u||u|) \biggr)
\biggl(-\frac{1}{\curpi} \biggr)\operatorname{Im}(a).\nonumber
\end{eqnarray}
By combining (\ref{e:IFT_spec_filter_d=0_log_on_the_offdiag}) and
(\ref{e:IFT_spec_filter_d=0}), one obtains the time domain kernels
on the right-hand side of~(\ref{e:time
domain_D=Jordan_and_eig=0}).

\section{Nonexistence of OFBM for certain exponents}

Proposition \ref{p:no_OFBM_with_H=Jordan_h=1}  is mentioned
in Remark \ref{rem:no_OFBM_with_H=Jordan_h=1}.

\begin{proposition}\label{p:no_OFBM_with_H=Jordan_h=1}
There does not exist an OFBM with exponent
\[
H = \pmatrix{ 1 & 0\cr
1 & 1
}.
\]
\end{proposition}
\begin{pf}
Assume that such an OFBM exists. For notational simplicity, denote
the process by $X$, and its entrywise processes by $X_1$ and
$X_2$. We will show that $X$ is not a proper process.

Note that, for $c > 0$, from the matrix expression for $c^{H}$ and o.s.s.,
%
\begin{equation}\label{e:o.s.s._root=1}
\left\{ \pmatrix{
X_{1}(ct)\cr
X_{2}(ct)
} \right\}_{t \in\bbR}
\stackrel{{\mathcal L}}= \left\{ \pmatrix{
c X_{1}(t) \cr
c \log(c) X_{1}(t) + c X_{2}(t)
} \right\}_{t \in\bbR}.
\end{equation}
In particular, this implies that $X_1$ is FBM with Hurst exponent
1. Thus, $X_1(t) = t Z$ a.s., where $Z$ is a Gaussian random
variable (e.g.,  \cite{taqqu2003}). By plugging this into
(\ref{e:o.s.s._root=1}),
we obtain
%
\begin{equation}\label{e:o.s.s._root=1_version2}
\left\{ \pmatrix{
ctZ \cr
X_{2}(ct)}
\right\}_{t \in\bbR}
\stackrel{{\mathcal L}}= \left\{\pmatrix{
ctZ \cr
c \log(c) tZ + c X_{2}(t)
} \right\}_{t \in\bbR}.
\end{equation}
In particular, by taking $c=t$ and $t=1$, (\ref{e:o.s.s._root=1_version2}) implies that
\[
EX_2(t)Z = E\bigl(t \log(t) Z + t X_2(1)\bigr)Z = t \log(t) EZ^2 + t EX_2(1)Z.
\]
Thus,
%
\begin{eqnarray} \label{e:o.s.s._root=1_increm}
&&E\bigl(X_2(t+h)-X_2(h)\bigr)\bigl(X_1(1+h)-X_1(h)\bigr)\nonumber\\
&&\quad = EX_2(t+h)Z - EX_2(h)Z
\\
&&\quad = \bigl((t+h) \log(t+h) - h \log(h)\bigr) EZ^2 + t EX_2(1)Z.\nonumber
\end{eqnarray}
By stationarity of the increments, the expression
(\ref{e:o.s.s._root=1_increm}) does not depend on $h$, which is
possible only when
\[
EZ^2 = 0.
\]
As a consequence, $\{X_2(ct)\}_{t \in\bbR} \stackrel{{\mathcal L}}= \{
c X_2(t)\}_{t \in\bbR}$. Thus,
\[
\pmatrix{
X_{1}(t)\cr
X_{2}(t)
}
= \pmatrix{
0 \cr
tY
}\qquad  \mbox{a.s.},
\]
where $Y$ is a Gaussian random variable. In particular, $X$ is not proper.
\end{pf}

\section{The exponential of a matrix in Jordan~canonical form} \label{a:Jordan_form}

Initially, let $J_{\lambda} \in M(n,\bbC)$ be a Jordan block of
size $n_\lambda$, whose expression is
%
\begin{equation} \label{e:Jordan_block}
J_{\lambda} = \pmatrix{
\lambda& 0 & 0 & \ldots& 0 \cr
1 & \lambda& 0 & \ldots& 0 \cr
0 & 1 & \lambda& \ldots& 0 \cr
\vdots& \vdots& \vdots& \ddots& \vdots\cr
0 & 0 & \ldots& 1 & \lambda
}.
\end{equation}
We have
%
\begin{equation} \label{e:z^Jlambda}
z^{J_{\lambda}} = \pmatrix{
z^{\lambda} & 0 & 0 & \ldots& 0\cr
(\log z)z^{\lambda} & z^{\lambda} & 0 & \ldots& 0\cr
\displaystyle\frac{(\log z)^{2}}{2!} z^{\lambda} & (\log z)z^{\lambda} & z^{\lambda}
& \ddots& 0 \cr
\vdots & \vdots & \ddots & \ddots & 0 \cr
\displaystyle\frac{(\log z)^{n_{\lambda}-1}}{(n_{\lambda}-1)!} z^{\lambda} &
\displaystyle\frac{(\log z)^{n_{\lambda}-2}}{(n_{\lambda}-2)!} z^{\lambda} &
\ldots& (\log z)z^{\lambda} & z^{\lambda}}.
\end{equation}
The expression for $z^{J}$, where $J$ is, more generally, a matrix
in Jordan canonical form (i.e., whose diagonal is made up of
Jordan blocks), follows immediately. In particular, the series-based
notion of the matrix exponential is consistent with the primary matrix
function-based notion of the matrix exponential.
\end{appendix}

\section*{Acknowledgments}

The first author was supported in part by the Louisiana Board of
Regents award LEQSF(2008-11)-RD-A-23.
The second author was supported in part by the NSF Grants
DMS-0505628 and DMS-0608669. The authors would like to thank
Professors\ Eric Renault and Murad Taqqu for their comments on this
work and to thank the two anonymous reviewers for their comments and
suggestions.

\printhistory


\begin{thebibliography}{00}

\bibitem{bahadoranbenassidebicki2003}
Bahadoran, C., Benassi, A. and D\c{e}bicki, K. (2003). Operator-self-similar
{G}aussian processes with stationary increments.
Preprint. Available at
\texttt{\href{http://math.univ-bpclermont.fr/prepublications/2003/2003-03.ps}{http://math.univ-bpclermont.fr/}\break
\href{http://math.univ-bpclermont.fr/prepublications/2003/2003-03.ps}{prepublications/2003/2003-03.ps}}.

\bibitem{beckerkernpap2008}
Becker-Kern, P. and Pap, G. (2008). Parameter estimation of selfsimilarity
exponents. \textit{J. Multivariate Anal.} \textbf{99} 117--140.
\MR{2408446}

\bibitem{biermemeerschaertscheffler2007}
Bierm{\'e}, H., Meerschaert, M.M. and Scheffler, H.-P. (2007). Operator
scaling stable random fields. \textit{Stochastic Process. Appl.}
\textbf{117} 312--332.
\MR{2290879}

\bibitem{chung2002}
Chung, C.-F. (2002). Sample means, sample autocovariances, and linear
regression of stationary multivariate long memory processes. \textit
{Econometric Theory} \textbf{18} 51--78.
\MR{1885350}

\bibitem{davidsondejong2000}
Davidson, J. and de Jong, R. (2000). The functional {C}entral {L}imit
{T}heorem and weak convergence to stochastic integrals {II}. \textit
{Econometric Theory} \textbf{16} 643--666.
\MR{1802836}

\bibitem{davidsonhashimzade2007}
Davidson, J. and Hashimzade, N. (2008). Alternative frequency and time domain
versions of fractional {B}rownian motion. \textit{Econometric Theory}
\textbf{24} 256--293.
\MR{2397267}

\bibitem{delgado2007}
Delgado, R. (2007). A reflected f{B}m limit for fluid models with {ON}/{OFF}
sources under heavy traffic. \textit{Stochastic Process. Appl.} \textbf{117} 188--201.
\MR{2290192}

\bibitem{didierpipiras2010}
Didier, G. and Pipiras, V. (2010). Exponents, symmetry groups and
classification of operator fractional {B}rownian motions. Preprint.

\bibitem{doladomarmol2004}
Dolado, J. and Marmol, F. (2004). Asymptotic inference results for
multivariate long-memory processes. \textit{Econom. J.} \textbf{7} 168--190.
\MR{2076631}

\bibitem{doob1953}
Doob, J.L. (1953). \textit{Stochastic Processes}. New York: Wiley.
\MR{0058896}

\bibitem{doukhan2003}
Doukhan, P., Oppenheim, G. and Taqqu, M. (2003). \textit{Theory and Applications
of Long-Range Dependence}. Boston, MA: Birkh\"{a}user.
\MR{1956041}

\bibitem{embrechtsmaejima2002}
Embrechts, P. and Maejima, M. (2002). \textit{Selfsimilar Processes}.
Princeton, NJ: Princeton Univ. Press.
\MR{1920153}

\bibitem{evardjafari1992}
Evard, J.-C. and Jafari, F. (1992). A complex {R}olle's theorem. \textit
{Amer. Math. Monthly} \textbf{99} 858--861.
\MR{1191706}

\bibitem{gradshteynryzhik2007}
Gradshteyn, I. and Ryzhik, I. (2007). \textit{Table of Integrals,
Series, and
Products}, 7th ed. New York, NY: Academic Press.
\MR{2360010}

\bibitem{hannan1970}
Hannan, E. (1970). \textit{Multiple Time Series}. New York: Wiley.
\MR{0279952}

\bibitem{helgasonpipirasabry2010}
Helgason, H., Pipiras, V. and Abry, P. (2010). Fast and exact synthesis of
stationary multivariate {G}aussian time series using circulant matrix
embedding. Preprint.

\bibitem{hornjohnson1991}
Horn, R. and Johnson, C. (1991). \textit{Topics in Matrix Analysis}.
New York, NY: Cambridge Univ. Press.
\MR{1091716}

\bibitem{hudsonmason1981JrMultAnalysis}
Hudson, W. and Mason, J. (1981). Operator-stable laws. \textit{J.
Multivariate Anal.} \textbf{11} 434--447.
\MR{0629799}

\bibitem{hudsonmason1982}
Hudson, W. and Mason, J. (1982). Operator-self-similar processes in a
finite-dimensional space. \textit{Trans. Amer. Math. Soc.} \textbf
{273} 281--297.
\MR{0664042}

\bibitem{jurekmason1993}
Jurek, Z. and Mason, J. (1993). \textit{Operator-Limit Distributions in
Probability Theory}. New York, NY: Wiley.
\MR{1243181}

\bibitem{konstantopouloslin1996}
Konstantopoulos, T. and Lin, S.J. (1996). Fractional {B}rownian approximations
of queueing networks. In \textit{Stochastic Networks. Lecture Notes in
Statistics} \textbf{117} 257--273. New York: Springer.
\MR{1466791}

\bibitem{laharohatgi1981}
Laha, R.G. and Rohatgi, V.K. (1981). Operator self-similar stochastic
processes in {$\mathbf{R}^{d}$}. \textit{Stochastic Process.
Appl.} \textbf{12} 73--84.
\MR{0632393}

\bibitem{maejimamason1994}
Maejima, M. and Mason, J. (1994). Operator-self-similar stable processes.
\textit{Stochastic Process. Appl.} \textbf{54} 139--163.
\MR{1302699}

\bibitem{majewski2003}
Majewski, K. (2003). Large deviations for multi-dimensional reflected
fractional {B}rownian motion. \textit{Stochastics} \textbf{75} 233--257.
\MR{1994908}

\bibitem{majewski2005}
Majewski, K. (2005). Fractional {B}rownian heavy traffic approximations of
multiclass feedforward queueing networks. \textit{Queueing Systems}
\textbf{50} 199--230.
\MR{2143121}

\bibitem{marinuccirobinson2000}
Marinucci, D. and Robinson, P. (2000). Weak convergence of multivariate
fractional processes. \textit{Stochastic Process. Appl.} \textbf{86} 103--120.
\MR{1741198}

\bibitem{masonxiao2002}
Mason, J. and Xiao, M. (2002). Sample path properties of
operator-self-similiar {G}aussian random fields. \textit{Theory
Probab. Appl.} \textbf{46} 58--78.
\MR{1968707}

\bibitem{meerschaertscheffler2001}
Meerschaert, M.M. and Scheffler, H.-P. (2001). \textit{Limit
Distributions for
Sums of Independent Random Vectors: Heavy Tails in Theory and
Practice}. New York: Wiley.
\MR{1840531}

\bibitem{pipirastaqqu2003}
Pipiras, V. and Taqqu, M. (2003). Fractional calculus and its
connections to
fractional {B}rownian motion. In \textit{Long-Range Dependence: Theory and
Applications} (P. Doukhan, G. Oppenheim and M.S. Taqqu, eds.)
165--201. Boston, MA: Birkh\"{a}user.
\MR{1956050}

\bibitem{robinson2008}
Robinson, P. (2008). Multiple local {W}hittle estimation in stationary
systems. \textit{Ann. Statist.} \textbf{36} 2508--2530.
\MR{2458196}

\bibitem{samorodnitskytaqqu1994}
Samorodnitsky, G. and Taqqu, M. (1994). \textit{Stable Non-Gaussian Processes:
Stochastic Models with Infinite Variance}. New York: Chapman \& Hall.
\MR{1280932}

\bibitem{stoevtaqqu2007}
Stoev, S. and Taqqu, M. (2006). How rich is the class of multifractional
{B}rownian motions. \textit{Stochastic Process. Appl.} \textbf{116} 200--221.
\MR{2197974}

\bibitem{taqqu2003}
Taqqu, M.S. (2003). Fractional {B}rownian motion and long range
dependence.
In \textit{Theory and Applications of Long-Range Dependence} (P.
Doukhan, G.
Oppenheim and M.S. Taqqu, eds.) 5--38. Boston: Birkh\"{a}user.
\MR{1956042}

\bibitem{yaglom1957}
Yaglom, A. (1957). Some classes of random fields in $n$-dimensional space,
related to stationary random processes. \textit{Theory Probab.
Appl.} \textbf{II} 273--320.
\MR{0094844}

\bibitem{yaglom1987}
Yaglom, A. (1987). \textit{Correlation {T}heory of {S}tationary and {R}elated
{R}andom {F}unctions. {V}olume I: {B}asic {R}esults}. New York: Springer.
\MR{0893393}

\end{thebibliography}
\end{document}